\documentclass[11pt,reqno]{amsart} 
\usepackage[margin=1in]{geometry}
\usepackage{amssymb}
\usepackage{mathtools}
\usepackage{mathrsfs}
\usepackage{enumitem}
\setlist{leftmargin=*,topsep=2pt,itemsep=2pt}

\newtheorem{theorem}{Theorem}[section]
\newtheorem{lemma}[theorem]{Lemma}
\newtheorem{proposition}[theorem]{Proposition}

\newtheorem{assumption}[theorem]{Assumption}
\newtheorem{remark}[theorem]{Remark}
\newtheorem{definition}[theorem]{Definition}

\numberwithin{equation}{section}

\DeclareMathOperator{\E}{\mathbb E}
\DeclareMathOperator{\Prob}{\mathbb P}

\DeclareMathOperator{\Cov}{\mathrm{Cov}}

\DeclareMathOperator{\dist}{dist}
\newcommand{\R}{\mathbb R}
\newcommand{\Z}{\mathbb Z}

\newcommand{\Om}{\Omega}
\newcommand{\eps}{\varepsilon}
\newcommand{\dd}{\mathrm d}

\newcommand{\wh}[1]{\widehat{#1}}
\let\le\leqslant
\let\ge\geqslant

\usepackage[hidelinks]{hyperref}
\usepackage{etoolbox}
\newcounter{Hequnique}
\AtBeginEnvironment{equation}{\stepcounter{Hequnique}}

\AtBeginEnvironment{subequations}{}
\AtEndEnvironment{subequations}{}

\begin{document}

\title[Random Schr\"odinger homogenization]{Stochastic Homogenization of Schr\"odinger Equations with Exponentially Mixing Sub-Gaussian Potentials}

\author{Fanfei Meng}
\address{College of Science, China Agricultural University, Beijing, China}
\email{meng\_fanfei@cau.edu.cn}
\thanks{F.~Meng was supported by the National Natural Science Foundation of China (Grant No.~12601455) and partially supported by the National Natural Science Foundation of China (Grant No.~12671284).}

\date{\today}



\begin{abstract}
We study $i\partial_t u_\eps=H_\eps u_\eps$, $H_\eps=-\Delta+V\left(\frac{x}{\eps},\omega\right)$, where $V$ is a stationary ergodic sub-Gaussian field with exponential mixing. 
We establish the optimal algebraic error rate $\eps^{\frac d2}$ in $L^2(\Omega;L^2_x)$ for every $d\ge1$, prove it sharp by a matching lower bound under a non degeneracy condition, and upgrade the mean square convergence of Zhang--Bal \cite{ZB2014} to almost sure convergence. 
The results are non perturbative in the amplitude of $V$ and require neither Gaussian marginals nor weak randomness. 
We further establish a Gaussian fluctuation central limit theorem with explicit covariance. 
The limiting fluctuation is the solution to the linear stochastic Schr\"odinger equation with additive spatial white noise 
\begin{equation*}
i\partial_t\mathcal G=H_0\mathcal G+\sqrt{\widehat R(0)}\,\bar u\,\xi,
\qquad H_0=-\Delta+\E[V],
\end{equation*}
where $\xi$ is a time frozen spatial white noise on $\R^d$ and $\mathcal{G}$ takes values in $H^{-\mathfrak{s}}(\R^d)$ for every $\mathfrak{s}>\frac{d}2$. 
The proofs combine Duhamel iteration with multiscale cluster expansion for the cumulants of smoothed linear functionals of the potential. 
\end{abstract}

\maketitle

\section{Introduction}

Let $(\Omega,\mathcal F,\Prob)$ be a probability space, and let $V(y,\omega)$ be a real valued stationary ergodic random field on $\R^d$ with mean
\begin{equation*}
\bar V:=\E[V(0,\omega)].
\end{equation*}
We study the random Schr\"odinger equation with rapidly oscillating stationary random potential. 
Define the microscopic and macroscopic Schr\"odinger operators
\begin{equation*}
H_\eps:=-\Delta+V\!\left(\frac{x}{\eps},\omega\right),\qquad H_0:=-\Delta+\bar V.
\end{equation*}
The microscopic equation is
\begin{equation}\tag{LSE$_\eps$}\label{eq:micro}
\begin{cases}
i\partial_t u_\eps=H_\eps u_\eps,\\
u_\eps(0,x)=u_0(x),
\end{cases}
\qquad\qquad (t,x)\in\R\times\R^d,
\end{equation}
with a real potential $V$ and initial data $u_0\in L^2(\R^d)$. 
As $\eps\to0$, the microscopic model \eqref{eq:micro} formally converges to the deterministic effective equation
\begin{equation}\tag{LSE$_0$}\label{eq:macro}
\begin{cases}
i\partial_t\bar u=H_0\bar u,\\
\bar u(0,x)=u_0(x),
\end{cases}
\qquad\qquad (t,x)\in\R\times\R^d.
\end{equation}
Define the fluctuation field
\begin{equation*}
w_\eps:=u_\eps-\bar u.
\end{equation*}
The core mathematical problems addressed in this work are the decay of $w_\eps$ as $\eps\to0$, the explicit convergence rate, and the statistical limiting behavior of the rescaled fluctuation $\eps^{-\frac d2}w_\eps$.

The stochastic fluctuation behavior is entirely governed by the covariance decay of the centered field
\begin{equation*}
a(y):=V(y)-\bar V,
\end{equation*}
whose covariance is
\begin{equation}\label{eq:cov}
R(y):=\E[a(y)a(0)],\qquad \widehat R(\xi)=\int_{\R^d}e^{-i\xi\cdot y}R(y)\dd y.
\end{equation}
Under exponential mixing, $R$ decays exponentially and $\widehat R$ is continuous. 
The critical parameter is $\widehat R(0)=\int_{\R^d}R(y)\dd y$, the variance of the spatial block average.

The main results of this paper are threefold. 
First, we prove the optimal algebraic rate
\begin{equation*}
\sup_{0\le t\le T}\|w_\eps(t)\|_{L^2(\Omega;L^2(\R^d))}\le C_{T,u_0}\eps^{\frac d2}.
\end{equation*}
Second, we prove almost sure convergence,
\begin{equation*}
\lim_{\eps\to0}\sup_{0\le t\le T}\|w_\eps(t)\|_{L^2(\R^d)}=0\qquad\text{a.s.}
\end{equation*}
Third, we prove a fluctuation central limit theorem, namely that the rescaled field $\eps^{-\frac d2}w_\eps$ converges in distribution to the Gaussian solution $\mathcal G$ of the linear random Schr\"odinger equation
\begin{equation}\tag{SPDE}\label{eq:spde}
\begin{cases}
i\partial_t\mathcal G=H_0\mathcal G+\sqrt{\widehat R(0)}\,\bar u(t,x)\,\xi(x),\\
\mathcal G(0,x)=0,
\end{cases}
\qquad\qquad (t,x)\in\R\times\R^d.
\end{equation}
Here $\xi$ is a time frozen spatial white noise on $\R^d$. 
The additive white noise enters because the fluctuation is driven, at the first Duhamel order, by the centered potential $a\left(\frac{x}{\eps}\right)$, see \eqref{eq:w-split}.

Throughout this paper, we assume that $V$ has almost surely continuous sample paths on $\R^d$, hence is bounded on compact sets, and that it satisfies Assumption \ref{ass:main} below. 
The sub-Gaussian tail in Assumption \ref{ass:main}(i) ensures that, for every $\eps>0$, the operator $H_\eps=-\Delta+V(\frac{x}{\eps},\omega)$ is essentially self-adjoint on $C_c^\infty(\R^d)$ by the Faris--Lavine theorem \cite[Theorem X.38]{RS1975}, so \eqref{eq:micro} is globally well-posed, see Remark \ref{rem:wellposed} for more details. 

The three results hold in every dimension $d\ge1$, with constants depending only on $d$ and uniform in $\bar V$; only the argument for the upper bound for the second moment depends on the dimension.
The transparent Fourier computation of this bound is carried out first for $d=1,2,3$, where the far region dyadic shell weight $2^{k(d-4)}$, which comes from Lemma \ref{lem:far} and Proposition \ref{prop:second}, converges geometrically precisely when $d<4$. 
On the far region the quadratic phase dominates and the time factor decays as $16\eps^4|\xi'|^{-4}$; multiplying by the shell volume $r^{d-1}\dd r$ gives the radial power $r^{d-5}$, so the tail $\int_\eps^\infty r^{d-5}\dd r$ fails to converge precisely at $d\ge4$.
The dyadic sum uses only the $L^\infty$ bound on $\widehat R$ rather than its smoothness, so exponential mixing does not remove this threshold. 
For $d\ge4$ the same bound is obtained by a physical space oscillatory integral argument, with $M=\lceil \frac{d}4\rceil+1$ time integrations by parts, developed in Appendix \ref{app:highdim}. 
The matching lower bound, the all moment estimate, the almost sure convergence and the fluctuation central limit theorem (CLT) use only the near region $|\xi'|\le\eps$, which is pure volume and is dimension independent. 

\subsection{Background and setting} 

Quantitative stochastic homogenization for divergence form elliptic operators traces back to the foundational work of Papanicolaou--Varadhan \cite{PapanicolaouVaradhan1979}, who laid the framework for partial differential equations with stationary ergodic random coefficients. 
Substantial progress was later achieved by Gloria--Neukamm--Otto in a quantitative program, which yields sharp variance estimates and optimal $L^2$ homogenization rates under spectral gap assumptions \cite{GNOInventiones2015,GloriaOtto2017}. 
The spectral gap condition encodes a strong decorrelation mechanism, and the quantitative ergodicity analysis under it, the corresponding large scale regularity theory, and the variance bounds available when correlations decay more slowly are developed in \cite{GNOInventiones2015,GNOregularity2020,GNOcorrelated2019}.
A separate perturbative harmonic-analytic line, initiated by Bourgain for weakly random discrete elliptic equations and extended to continuum mixing fields by Duerinckx--Lemm--Pagano \cite{DLP2025}, targets the enhanced accuracy of ensemble averages rather than pathwise fluctuations.

For parabolic and transport equations in random media, homogenization and the associated fluctuation central limit theorem have been extensively investigated by Bal \cite{Bal2010,Bal2011}, Gu--Bal \cite{GB2015}, and Bal--Jing \cite{BJ2010}, and related results for scalar homogenization problems and for the random heat equation appear in \cite{DunlapGuRyzhikZeitouni2021,Gu2016}. 
In this setting, the heat semigroup possesses intrinsic smoothing, which damps high frequency random oscillations automatically and makes the fluctuation argument comparatively tractable. 

Homogenization for Schr\"odinger type systems begins with the periodic setting, where the deterministic periodic theory and the effective mass theorems originate from the Bloch wave analysis of Allaire--Piatnitski \cite{AP2005}, and operator norm corrector estimates for periodic elliptic operators are developed in Birman--Suslina \cite{BirmanSuslina2006}. 
For stationary random potentials the primary prior quantitative contribution is due to Zhang--Bal \cite{ZB2014}, who established mean square convergence for large Gaussian random potentials, and further results for time dependent random Schr\"odinger homogenization appear in Gu--Ryzhik \cite{GuRyzhik2016}. 
Recent work of Harrop--Griffiths--Killip--Visan \cite{HKV2025} and Harrop--Griffiths--Ntekoume \cite{HGN2025} treats the nonlinear Schr\"odinger equation with sprinkled nonlinearity,
\begin{equation*}
i\partial_t\psi=-\partial_x^2\psi+2|\psi|^2\psi\dd\mu,\qquad x\in\R,
\end{equation*}
where $\mu$ is a homogeneous Poisson process. \cite{HKV2025} establishes global well-posedness, \cite{HGN2025} proves homogenization of this flow together with the fluctuations about the homogenized solution, and the choice $\mu=\delta_0$, which concentrates the nonlinearity at a single point, is treated in \cite{HKV2025b}. 

The unitary structure of the Schr\"odinger propagator changes the fluctuation computation.
Unlike the elliptic and parabolic equations, whose solution operators smooth away high frequencies, the Schr\"odinger group does not damp the off diagonal frequencies, and the fluctuations are not organized through a cell problem with an $L^2$ corrector.
The time double integral in the second moment is
\begin{equation*}
\int_0^t\int_0^t e^{i(r-r')\sigma}\dd r\dd r'
=\frac{2(1-\cos(t\sigma))}{\sigma^2},
\qquad
\sigma=|\xi|^2+|\eta|^2-2\xi\cdot\eta+\bar V=|\xi-\eta|^2+\bar V.
\end{equation*}
The time factor $\frac{2(1-\cos(t\sigma))}{\sigma^2}=\left|\int_0^t e^{ir\sigma}\dd r\right|^2$ in \eqref{eq:time-factor} is nonnegative and bounded, including at $\sigma=0$ where it equals $t^2$. Its oscillation in $\sigma$, rather than any positivity of the potential, is what localizes the spatial integral to the diagonal $\xi'=0$ and closes the estimate.
Because the propagator is unitary and the effective potential is just $\bar V=\E[V]$, no $L^2$ corrector enters the expansion, and the first Duhamel term already carries the full fluctuation while the remainder is of higher order. 
The limit is therefore a linear Schr\"odinger equation with additive spatial white noise, whose solution is the spatial white noise integral studied in \cite{DaPratoZabczyk2014,deBouardDebussche1999,OhPocovnicuWang2020}.
The quantitative objective is thus the same as that pursued by Gloria--Neukamm--Otto \cite{GNOInventiones2015,GloriaOtto2017,GNOregularity2020}, namely sharp variance estimates and optimal rates, but the machinery is not. Where that program builds on an elliptic corrector, a resolvent, and a spectral gap, the present argument relies on the oscillatory localization of the unitary group together with exponential mixing and sub-Gaussian tails.

We recall the standard definitions. 
Let $\tau_z$ denote the shift on $\Omega$, defined by
\begin{equation*}
V(\tau_z\omega,y):=V(\omega,y+z).
\end{equation*}

\begin{definition}[Stationarity]
The field $V$ is \emph{stationary} if the shift $\tau_z$ preserves $\Prob$ for every $z\in\R^d$, that is, for any $k\ge1$, any $y_1,\dots,y_k\in\R^d$, and any Borel set $A\subset\R^k$,
\begin{equation*}
\Prob\bigl((V(y_1),\dots,V(y_k))\in A\bigr)
=
\Prob\bigl((V(y_1+z),\dots,V(y_k+z))\in A\bigr)
\qquad\forall z\in\R^d.
\end{equation*}
\end{definition}

\begin{definition}[Ergodicity]
The stationary field $V$ is \emph{ergodic} if the shift group $\{\tau_z\}_{z\in\R^d}$ is ergodic, that is, any event $B\in\mathcal F$ satisfying $\tau_z^{-1}B=B$ for all $z\in\R^d$ has $\Prob(B)\in\{0,1\}$.
\end{definition}

To our knowledge, the only prior homogenization result for the dispersive Schr\"odinger equation with stationary random potentials is due to Zhang--Bal \cite{ZB2014}. 
Their analysis is specific to Gaussian potentials, for which the Duhamel expansion can be summed by a multiple scattering decomposition. 
Define $\|f\|_{d_1,d_2}:=\|\langle x\rangle^{d_1}\langle\nabla\rangle^{d_2}f\|_{L^2(\R^d)}$, where $\langle\nabla\rangle$ is the Fourier multiplier $(1+|\xi|^2)^{\frac12}$. 

\begin{theorem}[Zhang--Bal, {\cite[Theorem 1.1]{ZB2014}}]\label{th:zb}
Let $\mathfrak m\ge2$ and $d>\mathfrak m$. Let $q$ be a real valued, mean zero, stationary Gaussian field with correlation $R(x)=\E[q(y)q(x+y)]$ whose power spectrum $\widehat R$ is smooth and rapidly decaying, and assume $\|R\|_{12d,12d}<\infty$. 
Then, for every $T>0$, every $t\in(0,T)$, and every $u_0$ with $\langle\xi\rangle^{6d}\widehat u_0\in L^2(\R^d)$, the solution $u_\eps$ of $i\partial_tu_\eps=\big(P(D)-\eps^{-\frac{\mathfrak m}2}q(\frac{x}\eps)\big)u_\eps$ satisfies
\begin{equation*}
\lim_{\eps\to0}\E\Big[\|u_\eps(t)-\bar u(t)\|_{L^2(\R^d)}^2\Big]=0,
\end{equation*}
where $\bar u$ solves $i\partial_t\bar u=(P(D)-\rho)\bar u$ with $\rho=\int_{\R^d}\widehat R(\xi)\,|\xi|^{-\mathfrak m}\,d\xi$.
\end{theorem}

In this paper, we go beyond the mean square convergence of Theorem \ref{th:zb}, and we do so non perturbatively for potentials of every amplitude. 
We establish, to our knowledge for the first time for a dispersive Schrödinger equation, the complete quantitative picture, namely the sharp algebraic rate with a matching lower bound, the almost sure strong convergence, and a fluctuation central limit theorem with explicit covariance, under exponential mixing plus sub-Gaussian tails.
In particular, to our knowledge, the fluctuation limit is the first central limit theorem for the Schr\"odinger equation with a large stationary random potential, while the spatial white noise limit itself parallels the parabolic result of \cite{GB2015}.
Neither a Gaussian law nor a weakly random small amplitude assumption is needed.

We now state the standing assumption.

\begin{assumption}\label{ass:main}
The centered stationary ergodic random field $a(y,\omega):=V(y,\omega)-\bar V$, with $\bar V=\E[V(0)]$, satisfies the following two conditions.
\begin{enumerate}[label=(\roman*),font=\normalfont]
\item The marginals are sub-Gaussian, namely there exists $C_0>0$ such that
\begin{equation*}\E\bigl[e^{\lambda a(y)}\bigr]\le e^{\frac12 C_0\lambda^2}\qquad\forall\,y\in\R^d,\ \lambda\in\R.
\end{equation*}
\item The strong mixing ($\alpha$-mixing) coefficient decays exponentially, namely there exist $\kappa,\delta>0$ such that for all $R>0$,
\begin{equation}\label{eq:mix}
\bigl|\Prob(A\cap B)-\Prob(A)\Prob(B)\bigr|\le\delta e^{-\kappa R},
\end{equation}
whenever $A\in\sigma(a(y):|y|\le R)$ and $B\in\sigma(a(y):|y|\ge2R)$.
\end{enumerate}
\end{assumption}

Condition (i) gives uniform exponential integrability and closes the cumulant expansion, while condition (ii) is the exponential decay of the strong mixing ($\alpha$-mixing) coefficient and controls the off-diagonal terms.
The concentric balls in \eqref{eq:mix} are a convenient form of the standard $\alpha$-mixing coefficient. By stationarity the two formulations are equivalent up to a change of constants, since the defining $\sigma$-fields may be translated and enclosed by such balls.
Neither implies the other. 
Assumption \ref{ass:main} and the Gaussian class treated by Zhang--Bal \cite{ZB2014} are independent requirements, in that exponential mixing and a smooth, rapidly decaying spectrum neither implies nor is implied by the other.

The exponential rate in condition (ii) is used throughout. Under a merely polynomial mixing decay, the far region geometric sum is replaced by a summability condition on the mixing coefficient, the fluctuation scale and the central limit regime may change, and the present non perturbative estimates are not expected to close without further structure. These questions are left open.

\begin{remark}\label{rem:wellposed} 
Under Assumption \ref{ass:main}(i), the centered field admits almost surely a decomposition $a=B+a^{\sharp}$, by stationarity and the sub-Gaussian tail bound, where
\begin{equation*}
B(y):=a(y)\mathbf1_{\{|a(y)|\le1+|y|^2\}},\qquad
a^{\sharp}(y):=a(y)\mathbf1_{\{|a(y)|>1+|y|^2\}}.
\end{equation*}
By construction $|B(y)|\le1+|y|^2$. It remains to show that $a^{\sharp}\in L^p(\R^d)$ for every $p>0$. The tail bound $\Prob(|a(y)|>t)\le2e^{-\frac{t^2}{2C_0}}$ and stationarity give
\[
\int_{\R^d}\Prob\bigl(|a(y)|>1+|y|^2\bigr)\dd y<+\infty,
\]
and Tonelli's theorem yields
\[
\E\bigl[\|a^{\sharp}\|_{L^p}^p\bigr]\le C\int_{\R^d}\E\bigl[|a|^p\mathbf1_{\{|a|>1+|y|^2\}}\bigr]\dd y<+\infty,
\]
so $a^{\sharp}\in L^p(\R^d)$ almost surely.

Fix $\eps>0$. The potential $\bar V+B(\frac{x}\eps)$ is locally square integrable and bounded below by a quadratic function, while $a^{\sharp}(\frac{x}\eps)\in L^p(\R^d)$ for $p>\max\{2,\frac{d}2\}$. The Faris--Lavine theorem \cite[Theorem X.38]{RS1975}, together with the $L^p$ relative boundedness of $a^{\sharp}(\frac{x}{\eps})$, then shows that $H_\eps=-\Delta+V(\frac{x}{\eps},\omega)$ is essentially self-adjoint on $C_c^\infty(\R^d)$. With $H_\eps$ denoting its self-adjoint closure, Stone's theorem gives the unitary group $e^{-itH_\eps}$, so \eqref{eq:micro} is globally well-posed in $L^2(\R^d)$ with $u_\eps(t)=e^{-itH_\eps}u_0$ and $\|u_\eps(t)\|_{L^2}=\|u_0\|_{L^2}$. The argument is standard, see \cite[Section X.3]{RS1975}.
\end{remark}

Four concrete examples illustrate the relationship between the two assumptions and why both are needed.

\begin{enumerate}
\item[(a)] Let $W$ be white noise on $\R^d$. 
The stationary Gaussian field
\begin{equation*}
	V_1(y):=\int_{\R^d}e^{iy\cdot\xi}(1+|\xi|^2)^{-\frac{d+\alpha}{4}}\,\dd W(\xi),\qquad \alpha>d,
\end{equation*}
has covariance decaying as $|y|^{-\alpha}$, which is polynomially rather than exponentially decaying. 
It is sub-Gaussian but fails exponential mixing, because exponential mixing forces exponentially decaying covariance, whereas the covariance here decays only polynomially. Since $\alpha>d$ the covariance is integrable and $\widehat R(0)<\infty$. 
This shows that sub-Gaussianity alone does not give quantitative homogenization. 

\item[(b)] Let $\{\eta_k\}_{k\in\Z^d}$ be i.i.d.\ Cauchy variables and $\chi\in C_c^\infty((-\frac12,\frac12)^d)$ a bump. 
The field
\begin{equation*}
V_2(y):=\sum_{k\in\Z^d}\eta_k\,\chi(y-k)
\end{equation*}
is perfectly mixing (disjoint supports) but has no exponential moment. 
This shows that mixing alone does not control higher cumulants.

\item[(c)] The Gaussian field
\begin{equation*}
V_3(y):=\int_{\R^d}e^{iy\cdot\xi}(1+|\xi|^2)^{-\frac12}\,\dd W(\xi)
\end{equation*}
has covariance given by the Bessel kernel, decaying as $e^{-c|y|}$. It satisfies both assumptions and is the canonical model covered by our results, with $\widehat R(0)=(2\pi)^d>0$. 
Its sample paths are H\"older continuous for $d=1$ and degenerate as $d$ grows. The regularity is borderline (logarithmic) at $d=2$ and distributional for $d\ge3$. For $d\ge2$ the standing continuity hypothesis should therefore be read as holding for a smooth approximation of $V_3$, while the main results hold for every $d\ge1$ under Assumption \ref{ass:main}.

\item[(d)] Let $\eta$ be Gaussian with exponentially decaying covariance and $f=\tanh$. 
The bounded non-Gaussian field
\begin{equation*}
V_4(y):=\tanh(\eta(y))
\end{equation*}
satisfies both assumptions by Lipschitz transfer of mixing 
\begin{equation*}
|\Cov(g(V_4),h(V_4))|\le\|g\|_\infty\|h\|_\infty e^{-\kappa R}.
\end{equation*}
This is the standard non-Gaussian sub-Gaussian example.
\end{enumerate}

The centered field is $a(y)=V(y)-\bar V$, and the covariance $R$ and the spectral density $\wh R$ are those defined in \eqref{eq:cov}.
Assumption \ref{ass:main} implies $R\in L^1(\R^d)$ and $\wh R(0)\ge0$, with $\wh R$ bounded and continuous.

\subsection{Methods and results}

The Schr\"odinger group generated by $H_0$ is
\begin{equation*}S(t,r):=e^{-iH_0(t-r)},\qquad \mathcal F(S(t,r)f)(\xi)=e^{-i(|\xi|^2+\bar V)(t-r)}\widehat f(\xi).
\end{equation*}
Subtracting \eqref{eq:macro} from \eqref{eq:micro}, the fluctuation equation is
\begin{equation*}i\partial_t w_\eps=H_0 w_\eps+a_\eps(\bar u+w_\eps),\qquad w_\eps(0)=0,
\end{equation*}
where $a_\eps(x):=a\!\left(\frac{x}{\eps}\right)$. 
The Duhamel formula shows that 
\begin{equation*}w_\eps(t)=-i\int_0^t S(t,r)a_\eps(r)\bigl(\bar u(r)+w_\eps(r)\bigr)\dd r.
\end{equation*}
We split $w_\eps=w_\eps^{(1)}+w_\eps^{(2)}$ by the first Duhamel iteration,
\begin{equation}\label{eq:w-split}
\begin{aligned}
w_\eps^{(1)}(t)&:=-i\int_0^t S(t,r)a_\eps(r)\bar u(r)\dd r,\\
w_\eps^{(2)}(t)&:=-i\int_0^t S(t,r)a_\eps(r)w_\eps(r)\dd r.
\end{aligned}
\end{equation}
The second term is implicit, since the full fluctuation $w_\eps$ occurs under its integral. Expanding that $w_\eps$ by the Duhamel formula once more and repeating the self-similar substitution resolves the equation into the hierarchy
\begin{equation*}
w_\eps=\sum_{k\ge1}w_\eps^{(k)},
\end{equation*}
where $w_\eps^{(k)}$ carries $k$ ordered insertions of the centered field, with $w_\eps^{(1)}$ the linear term and $r_\eps:=\sum_{k\ge2}w_\eps^{(k)}$ the remainder. The two parts are controlled by different estimates.

\begin{itemize}
\item For $w_\eps^{(1)}$, its Fourier transform is the spatial average of the generalized field $\widehat a(\eps\eta)$ against the free profile $\widehat u_0(\xi-\eta)$ and the time factor $\int_0^t e^{ir\sigma}\dd r$, with $\sigma=2\xi\cdot\eta-|\eta|^2$. Set $\zeta=\xi-\eta$ and then $\xi'=\eps\eta$, so that $\xi=\zeta+\frac{\xi'}{\eps}$ and $\dd\xi=\eps^{-d}\dd\xi'$. By Plancherel, the second moment is
\begin{equation*}
\E\left[\|w_\eps^{(1)}(t)\|_{L^2_x}^2\right]
=(2\pi)^{-2d}\int_{\R^d}|\widehat u_0(\zeta)|^2
\int_{\R^d}\widehat R(\xi')\frac{2(1-\cos(t\phi))}{\phi^2}\dd\xi'\dd\zeta,
\qquad
\phi=\frac{2\zeta\cdot\xi'}{\eps}+\frac{|\xi'|^2}{\eps^2}.
\end{equation*}
The $\xi'$ integral is split at $|\xi'|=\eps$. On the near region $|\xi'|\le\eps$, the phase is bounded and $\frac{2(1-\cos(t\phi))}{\phi^2}\le t^2$, over a volume $O(\eps^d)$. On the far region $|\xi'|>\eps$, the quadratic term $\frac{|\xi'|^2}{\eps^2}$ dominates and $\frac{2(1-\cos(t\phi))}{\phi^2}\le\frac{16\eps^4}{|\xi'|^4}$; the transition strip $\eps<|\xi'|<4R\eps$ is counted as the near region, and a dyadic shell decomposition with shell weight $2^{k(d-4)}$ converges geometrically for $d\le3$. Combining the two regions,
\begin{equation}\label{eq:core-intro}
\E\left[\|w_\eps^{(1)}(t)\|_{L^2_x}^2\right]\le C_{T,u_0}\eps^d.
\end{equation}
The same representation restricted to the near region identifies the diagonal leading term, yields the matching lower bound under $\widehat R(0)>0$ and $\widehat u_0(0)\ne0$, and fixes the fluctuation covariance.

\item For $w_\eps^{(k)}$, $k\ge2$, Lemma \ref{lem:iterate} gives
\begin{equation*}
\E\left[\|w_\eps^{(k)}(t)\|_{L^2_x}^2\right]\le C^k\eps^{dk},
\qquad k\ge2,
\end{equation*}
and the estimates are summable over $k$, the first term being the largest,
\begin{equation}\label{eq:rem-intro}
\|r_\eps(t)\|_{L^2(\Om;L^2_x)}^2
\le\left(\sum_{k\ge2}C^{\frac k2}\eps^{\frac{dk}{2}}\right)^2
\le C\eps^{2d}.
\end{equation}
The remainder is therefore of order $\eps^{d}$, one half power smaller than the leading $\eps^{\frac d2}$ term.
\end{itemize}

Combining the leading bound \eqref{eq:core-intro} with the smaller remainder \eqref{eq:rem-intro} by the triangle inequality, and summing the Duhamel series, yields the first main result, the optimal algebraic convergence rate \(\eps^{\frac d2}\) in \(L^2(\Om;L^2_x)\).

\begin{theorem}\label{th:rate}
Under Assumption \ref{ass:main}, the following hold for every $d\ge1$.
\begin{enumerate}[label=(\roman*),font=\normalfont]
\item For all $0<\eps\le1$,
\begin{equation*}
\sup_{0\le t\le T}\|u_\eps(t)-\bar u(t)\|_{L^2(\Om;L^2(\R^d))}\le C_{T,u_0}\eps^{\frac d2}.
\end{equation*}
If additionally $\widehat R(0)>0$ and $\widehat{u_0}(0)\neq0$, then there exist $c_{t,u_0}>0$ and $\eps_0>0$ such that
\begin{equation}\label{eq:lower-intro}
\E\left[\|u_\eps(t)-\bar u(t)\|_{L^2(\R^d)}^2\right]\ge c_{t,u_0}\,\eps^d
\qquad\text{for all }0<\eps\le\eps_0.
\end{equation}
The rate $\eps^{\frac d2}$ is sharp.

\item For every finite $m\ge1$,
\begin{equation*}
\E\left[\sup_{0\le t\le T}\|u_\eps(t)-\bar u(t)\|_{L^2(\R^d)}^m\right]\le C_{m,T,u_0}\eps^{\frac{md}{2}}.
\end{equation*}
\end{enumerate}
\end{theorem}

\begin{remark}
	The lower bound \eqref{eq:lower-intro} requires $\widehat R(0)>0$.
	If $\widehat R(0)=0$, the block-average variance is governed by the next spectral moment,
	the $\eps^{\frac{d}2}$ fluctuations are suppressed, and the rate and CLT scale must be recomputed from the local regularity of $\widehat R$ near the origin. 
	We do not treat this case.
	The point condition $\widehat u_0(0)\ne0$ is a concise way of requiring $\widehat u_0$ to be nonzero in a neighborhood of the origin, which is all that the lower bound proof uses; it is not essential that the nonzero value be attained exactly at the origin.
\end{remark}

\begin{remark}
	The compact support of $\widehat u_0$ is technical and localizes the $\zeta$-integral only.
	For $d\le3$, the constants depend only on $\|u_0\|_{L^2}$, so the estimates extend to $u_0\in L^2(\R^d)$ by density directly. 
	For $d\ge4$, the physical space bound of Lemma \ref{lem:fresnel} is an operator norm estimate uniform in $\eps$, and the $\eps^d$ prefactor absorbs the Schwartz norm dependence. 
	Approximating $u_0\in L^2$ by $C_c^\infty$ converges in $L^2(\Omega;L^2_x)$ by the unitarity of $S(t,r)$.
\end{remark} 

\begin{remark} 	
	The sub-Gaussian assumption \ref{ass:main}(i) is used at two levels. The $L^2$ rate requires only a moment of order $p>2$ together with exponential mixing, while the fluctuation CLT and the all moment estimate require the full sub-Gaussian proxy to close the cumulant expansion and apply Proposition \ref{prop:cluster-moment}.
\end{remark} 

The rate estimate above is in $L^2(\Om;L^2_x)$. 
To upgrade it to almost sure convergence, we apply the dyadic Borel--Cantelli argument. 
At the dyadic sequence $\eps_n=2^{-n}$,
\begin{equation*}
\E\left[\sup_{0\le t\le T}\|w_{\eps_n}(t)\|_{L^2_x}^{2q}\right]\le C^q\eps_n^{dq}.
\end{equation*}
By Markov's inequality,
\begin{equation*}
\Prob\left(\sup_t\|w_{\eps_n}(t)\|>\eps_n^{\frac d2-\eta}\right)\le C_q\eps_n^{2q\eta},
\end{equation*}
summable for $q$ large. 
The time increment satisfies
\begin{equation*}
\E\left[\|w_\eps^{(1)}(t)-w_\eps^{(1)}(s)\|_{L^2_x}^{2q}\right]\le C|t-s|^{2q},
\end{equation*}
so Kolmogorov--Chentsov gives equicontinuity in $C([0,T];L^2_x)$ uniformly in $\eps$, filling the gaps between dyadic points. 
This yields the almost sure strong convergence.

\begin{theorem}\label{th:as}
Under Assumption \ref{ass:main}, for every $d\ge1$, every fixed $T>0$ and $u_0\in L^2(\R^d)$,
\begin{equation*}
\lim_{\eps\to0}\sup_{0\le t\le T}\|u_\eps(t)-\bar u(t)\|_{L^2(\R^d)}=0
\qquad\text{a.s.}
\end{equation*}
\end{theorem}

Beyond the convergence, we study the limiting distribution of the rescaled fluctuation. 
Define
\begin{equation*}
\mathcal F_\eps(t,x):=\eps^{-\frac d2}w_\eps(t,x),
\qquad
W_t^\varphi(x):=\int_0^t \overline{(S(t,r)^*\varphi)(x)}\,\bar u(r,x)\,\dd r.
\end{equation*}
By \eqref{eq:w-split} and the remainder estimate \eqref{eq:rem-intro},
\begin{equation*}
\langle\mathcal F_\eps(t),\varphi\rangle
=-i\eps^{-\frac d2}\int_{\R^d}a\!\left(\frac{x}{\eps}\right)W_t^\varphi(x)\,\dd x+o_{L^2(\Om)}(1),
\end{equation*}
so the rescaled fluctuation is, to leading order, a weighted spatial average of the centered field. 
Its variance converges to the explicit limit
\begin{equation*}
\lim_{\eps\to0}\E\left[\langle\mathcal F_\eps(t),\varphi\rangle\overline{\langle\mathcal F_\eps(s),\psi\rangle}\right]
=\widehat R(0)\int_{\R^d}W_t^\varphi(x)\,\overline{W_s^\psi(x)}\,\dd x,
\end{equation*}
and the higher cumulants vanish under the fluctuation scaling by the cluster estimate \eqref{eq:cluster-moment}, as recorded in Remark \ref{rem:cluster-boundary}. 
The finite dimensional distributions are therefore Gaussian with the covariance \eqref{eq:cltcov}, and tightness in $C([0,T];H^{-\mathfrak s}(\R^d))$, $\mathfrak s>\frac d2$, follows from the second and fourth moment increment estimates of Proposition \ref{prop:tight}. 
This yields the final main result, the Gaussian fluctuation limit.

\begin{theorem}\label{th:clt}
Assume $d\ge1$ and fix a Sobolev exponent $\mathfrak s>\frac d2$.
As $\eps\to0$, the family $\{\mathcal F_\eps\}_{\eps>0}$ converges in distribution in $C([0,T];H^{-\mathfrak s}(\R^d))$ to a centered Gaussian field $\mathcal G$ with covariance
\begin{equation}\label{eq:cltcov}
\Cov\bigl(\langle\mathcal G(t),\varphi\rangle,\langle\mathcal G(s),\psi\rangle\bigr)
= \wh R(0)\int_{\R^d} W_t^\varphi(x)\,\overline{W_s^\psi(x)}\,\dd x,
\qquad \varphi,\psi\in H^{\mathfrak s}(\R^d),
\end{equation}
where
\begin{equation*}
W_t^\varphi(x):=\int_0^t \overline{(S(t,r)^*\varphi)(x)}\,\bar u(r,x)\,\dd r.
\end{equation*}
\end{theorem}

\begin{remark}
The covariance \eqref{eq:cltcov} is exactly that of the mild solution to \eqref{eq:spde}, namely
\[
\mathcal G(t)=-i\sqrt{\widehat R(0)}\int_0^t S(t,r)\bigl(\bar u(r)\,\xi\bigr)\dd r,
\]
where $\xi$ is a time frozen spatial white noise on $\R^d$, the centered Gaussian field with covariance $\E[\xi(x)\xi(y)]=\delta(x-y)$.
The process is a well-defined Gaussian variable in every negative Sobolev space $H^{-\mathfrak s}(\R^d)$ with $\mathfrak s>\frac d2$. Writing $J=(1-\Delta)^{\frac12}$, the regularization $J^{-\mathfrak s}$ makes the spatial integral Hilbert--Schmidt, and
\[
\E\bigl[\|\mathcal G(t)\|_{H^{-\mathfrak s}}^2\bigr]
=\widehat R(0)\int_0^t\!\!\int_0^t
\bigl\langle J^{-\mathfrak s}S(t,r)\bar u(r),\ J^{-\mathfrak s}S(t,r')\bar u(r')\bigr\rangle_{L^2_x}\,\dd r\,\dd r'<\infty,
\]
finite since $J^{-\mathfrak s}$ is Hilbert--Schmidt precisely when $2\mathfrak s>d$ and $S(t,r)$ is unitary. The Cameron--Martin space of the underlying cylindrical Gaussian measure is $L^2(\R^d)$, and almost every sample of $\mathcal G$ lies in $H^{-\mathfrak s}$.

The limit \eqref{eq:spde} is the dispersive counterpart of \cite[Theorem 2.4]{GB2015}, the same time frozen spatial white noise, with the parabolic profile $u_{\hom}$ replaced by the free Schr\"odinger profile $\bar u$.

The limit does not depend on the chosen exponent. For any $\mathfrak s_1,\mathfrak s_2>\frac d2$, the restrictions of the two limits to $H^{-\mathfrak s_1}\cap H^{-\mathfrak s_2}$ coincide, since they are the same cylindrical Gaussian field on the common test space $H^{\mathfrak s_1}\cap H^{\mathfrak s_2}$, with the covariance \eqref{eq:cltcov}.
\end{remark}

The rest of the paper is organized as follows. 
Section \ref{sec:prelim} collects the analytic conventions and the probabilistic tools used throughout.
Section \ref{sec:second} develops the Duhamel decomposition and iteration. Subsection \ref{sec:fourier} gives the second moment estimates, Subsection \ref{sec:physical} the second order correction and the sharp lower bound, and Subsection \ref{sec:moments} the all moment estimates, which together prove Theorem \ref{th:rate}. The Fourier near/far second moment bound in Subsection \ref{sec:fourier} is written for $d\le3$.
Section \ref{sec:as} proves the almost sure convergence, Theorem \ref{th:as}.
Section \ref{sec:clt} proves the fluctuation central limit theorem, Theorem \ref{th:clt}.
Appendix \ref{app:cluster} gives the multiscale cluster expansion for linear functionals, and Appendix \ref{app:highdim} supplies the physical space oscillatory integral argument that extends the second moment upper bound to every $d\ge4$, the remaining estimates being dimension independent.

\section{Preliminaries}\label{sec:prelim}

We write $\mathcal F(g)(\xi)$ for the Fourier transform of a compound expression $g$, while $\widehat g$ denotes the Fourier transform of a single symbol. 
We use the convention
\begin{equation*}\mathcal F(f)(\xi)=\int_{\R^d}e^{-ix\cdot\xi}f(x)\dd x,\qquad
f(x)=(2\pi)^{-d}\int_{\R^d}e^{ix\cdot\xi}\mathcal F(f)(\xi)\dd\xi,
\end{equation*}
so Plancherel identity is 
\begin{equation*}\|f\|_{L^2}^2=(2\pi)^{-d}\|\mathcal F(f)\|_{L^2}^2.
\end{equation*}

The Schr\"odinger group generated by $H_0=-\Delta+\bar V$ is
\begin{equation}\label{eq:flow}
\mathcal F(S(t,r)f)(\xi)=e^{-i(|\xi|^2+\bar V)(t-r)}\widehat f(\xi),
\end{equation}
and its kernel is the Fresnel kernel
\begin{equation*}
K_t(x)=(4\pi i t)^{-\frac d2}e^{i\frac{|x|^2}{4t}},\qquad t\neq0.
\end{equation*}
The Fresnel kernel is understood as a tempered distribution for every $d$, defined through the oscillatory integral and characterized by $\mathcal F(K_t)(\xi)=e^{-it|\xi|^2}$, rather than as an ordinary integrable function.
The Schr\"odinger group $S(t,r)=e^{-iH_0(t-r)}$ is unitary on $L^2(\R^d)$ \cite[Section X.8]{RS1975}, that is,
\begin{equation*}
\|S(t,r)f\|_{L^2_x}=\|f\|_{L^2_x}
\end{equation*}
for all $t,r\in\R$ and $f\in L^2(\R^d)$.

\begin{lemma}[L\'evy continuity theorem, {\cite[Section 3.3]{Durrett2019}}] \label{lem:levy}
Let $\{\mu_n\}$ be probability measures on $\R^k$ with characteristic functions $\phi_n(\theta)=\int_{\R^k}e^{i\theta\cdot x}\dd\mu_n(x)$. 
If $\phi_n(\theta)$ converges pointwise to a function $\phi(\theta)$ that is continuous at $\theta=0$, then $\phi$ is the characteristic function of a probability measure $\mu$ on $\R^k$, and $\mu_n$ converges weakly to $\mu$.
\end{lemma}

\begin{lemma}[Borel--Cantelli lemma, {\cite[Section 2.3]{Durrett2019}}] \label{lem:bc}
If
\begin{equation*}
\sum_{n=1}^{\infty} \Prob(A_n)<\infty,
\end{equation*}
then, 
\begin{equation*}
\Prob \left( \limsup_{n\to\infty} A_n \right)=0.
\end{equation*}
\end{lemma}

\begin{lemma}[Isserlis' theorem, {\cite[Chapter~1]{Janson1997}}] \label{lem:isserlis}
If $X_1,\dots,X_{2q}$ are jointly centered Gaussian random variables, then,
\begin{equation*}
\E\left[\prod_{j=1}^{2q}X_j\right]
=\sum_{\pi\in\mathcal P(2q)}\prod_{\{j,k\}\in\pi}\E[X_jX_k],
\end{equation*}
where $\mathcal P(2q)$ is the set of all pairings of $\{1,\dots,2q\}$.
\end{lemma}

\begin{lemma}[Kolmogorov--Chentsov compactness criterion, {\cite[Theorem~3.23]{Kallenberg2002}}] \label{lem:kc}
Let $\{X_\eps(t)\}_{\eps>0}$ be a family of processes taking values in a separable Banach space $B$. If there exist constants $p\ge1$, $\alpha>1$, and $C_{\mathrm{inc}}>0$ such that
\begin{equation*}
\E\left[\|X_\eps(t)-X_\eps(s)\|_B^p\right]\le C_{\mathrm{inc}}|t-s|^\alpha
\end{equation*}
uniformly in $\eps$, then for every $\beta$ with $0<\beta<\frac{\alpha-1}{p}$ each $X_\eps$ admits a $\beta$-Hölder continuous modification, and the family is relatively compact in $C([0,T];B)$ whenever $\{X_\eps(0)\}$ is relatively compact in $B$.
The criterion is commonly stated for finite dimensional $B$; its dyadic proof uses only the separability of $B$ together with the Arzelà--Ascoli theorem, so the Banach valued form follows without further assumptions.
\end{lemma}

\begin{lemma}[Prokhorov, {\cite[Chapter~5]{Billingsley1999}}]\label{lem:prokhorov}
A family $\Pi$ of probability measures on a complete separable metric space $E$ is relatively compact in the weak topology if and only if it is tight, that is, for every $\delta>0$ there exists a compact $K_\delta\subset E$ such that
\begin{equation*}
\sup_{\mu\in\Pi}\mu(E\setminus K_\delta)<\delta.
\end{equation*}
\end{lemma}

\section{Duhamel decomposition and iteration}\label{sec:second}

Write $u_\eps=\bar u+w_\eps$. 
Subtracting \eqref{eq:macro} from \eqref{eq:micro},
\begin{equation*}i\partial_t w_\eps=H_0 w_\eps+a_\eps(\bar u+w_\eps),\qquad w_\eps(0)=0,
\end{equation*}
where we recall $a_\eps(x)=a\left(\frac{x}{\eps}\right)$, so that
\begin{equation}\label{eq:acov}
\E[a_\eps(y)a_\eps(y')]=R\biggl(\frac{y-y'}{\eps}\biggr)
\end{equation}
by \eqref{eq:cov}.
The Duhamel form is
\begin{equation}\label{eq:dw}
w_\eps(t)=-i\int_0^t S(t,r)a_\eps(r)\bigl(\bar u(r)+w_\eps(r)\bigr)\dd r.
\end{equation}
Iterating once gives the decomposition displayed in \eqref{eq:w-split}, with $w_\eps^{(1)}$ the linear term and $w_\eps^{(2)}$ the remainder.

\subsection{Second moment estimates}\label{sec:fourier}

The Schr\"odinger propagator $e^{i(t-s)\Delta}$ has infinite propagation speed. Its Fresnel kernel is nonzero for every $x\neq y$ with only algebraic decay. 
A direct physical space Gronwall estimate would therefore accumulate contributions from the entire space at every time step. 
We instead work in Fourier space, where the off diagonal terms are controlled by the rapid oscillation of the phase factor $\exp\left\{i(t-s)\left(\left|\xi+\frac{\xi'}{\eps}\right|^2-|\xi|^2\right)\right\}$ combined with the exponential decay of the spectral density $\widehat R(\xi')$.

By \eqref{eq:flow} and the convolution theorem, the Fourier transform of the first order correction has the following representation.

\begin{proposition}\label{prop:fourier}
Under Assumption \ref{ass:main}, for the first order Duhamel term $w_\eps^{(1)}$,
\begin{equation}\label{eq:fourier-w}
\widehat w_\eps^{(1)}(t,\xi)
=-i\eps^d e^{-i(|\xi|^2+\bar V)t}(2\pi)^{-d}
\int_{\R^d}\widehat a(\eps\eta)\widehat u_0(\xi-\eta)
\left[\int_0^t e^{ir\sigma(\xi,\eta)}\dd r\right]\dd\eta,
\end{equation}
where
\begin{equation*}\sigma(\xi,\eta):=|\xi|^2-|\xi-\eta|^2=2\xi\cdot\eta-|\eta|^2.
\end{equation*}
\end{proposition} 

\begin{proof}
By \eqref{eq:flow} and the convolution theorem $\mathcal F(fg)=(2\pi)^{-d}\widehat f*\widehat g$,
\begin{equation*}\widehat w_\eps^{(1)}(t,\xi)
=-i\eps^d\int_0^t\int_{\R^d}
e^{-i(|\xi|^2+\bar V)(t-r)}
\widehat a(\eps\eta)\,\mathcal F(\bar u(r))(\xi-\eta)\,(2\pi)^{-d}\,\dd\eta\,\dd r.
\end{equation*}
Since $\bar u$ solves \eqref{eq:macro},
\begin{equation*}\mathcal F(\bar u(r))(\xi')=e^{-i(|\xi'|^2+\bar V)r}\widehat u_0(\xi').
\end{equation*}
Substituting and integrating the time exponential,
\begin{equation*}
\int_0^t e^{-i(|\xi|^2+\bar V)(t-r)}e^{-i(|\xi-\eta|^2+\bar V)r}\dd r
=e^{-i(|\xi|^2+\bar V)t}\int_0^t e^{ir(|\xi|^2-|\xi-\eta|^2)}\dd r.
\end{equation*}
This gives \eqref{eq:fourier-w}.
\end{proof}

For the generalized Fourier transform of the stationary field,
\begin{equation*}\E[\widehat a(\alpha)\overline{\widehat a(\beta)}]
=(2\pi)^d\delta(\alpha-\beta)\widehat R(\alpha),
\end{equation*}
and with $\alpha=\eps\eta$, $\beta=\eps\eta'$,
\begin{equation}\label{eq:covfourier2}
\E[\widehat a(\eps\eta)\overline{\widehat a(\eps\eta')}]
=(2\pi)^d\eps^{-d}\delta(\eta-\eta')\widehat R(\eps\eta).
\end{equation}

We now compute the second moment in Fourier space. 
Multiplying the representation \eqref{eq:fourier-w} of Proposition \ref{prop:fourier} by its conjugate and taking expectation,
\begin{equation*}
\E\left[|\widehat w_\eps^{(1)}(t,\xi)|^2\right]
=\eps^d(2\pi)^{-d}\int_{\R^d}\widehat R(\eps\eta)
\left|\int_0^t e^{ir\sigma(\xi,\eta)}\dd r\right|^2
|\widehat u_0(\xi-\eta)|^2\dd\eta.
\end{equation*}
By Plancherel,
\begin{equation*}
\E\left[\|w_\eps^{(1)}(t)\|_{L^2_x}^2\right]
=\eps^d(2\pi)^{-2d}\int_{\R^d}\int_{\R^d}
\widehat R(\eps\eta)
\left|\int_0^t e^{ir\sigma(\xi,\eta)}\dd r\right|^2
|\widehat u_0(\xi-\eta)|^2\dd\eta\dd\xi.
\end{equation*}
The time double integral evaluates to
\begin{equation}\label{eq:time-factor}
\left|\int_0^t e^{ir\sigma}\dd r\right|^2
=\int_0^t\!\!\int_0^t e^{i(r-r')\sigma}\dd r\dd r'
=\frac{2(1-\cos(t\sigma))}{\sigma^2}.
\end{equation}
Change variables $\zeta=\xi-\eta$, so that $\sigma(\xi,\xi-\zeta)=|\xi|^2-|\zeta|^2$, and then set $\xi'=\eps(\xi-\zeta)$, so $\xi=\zeta+\frac{\xi'}{\eps}$ and $\dd\xi=\eps^{-d}\dd\xi'$. Then
\begin{equation*}
|\xi|^2-|\zeta|^2=\frac{2\zeta\cdot\xi'}{\eps}+\frac{|\xi'|^2}{\eps^2},
\end{equation*}
and the integral becomes
\begin{equation}\label{eq:key-integral}
\E\left[\|w_\eps^{(1)}(t)\|_{L^2_x}^2\right]
=(2\pi)^{-2d}\int_{\R^d}|\widehat u_0(\zeta)|^2
\int_{\R^d}\widehat R(\xi')\,\frac{2(1-\cos(t\phi))}{\phi^2}\,\dd\xi'\,\dd\zeta,
\qquad
\phi:=\frac{2\zeta\cdot\xi'}{\eps}+\frac{|\xi'|^2}{\eps^2}.
\end{equation}
Since $\widehat u_0$ is compactly supported, write $\operatorname{supp}(\widehat u_0)\subset B_R$, so that $|\zeta|\le R$ uniformly. 
The estimate of the $\xi'$-integral in \eqref{eq:key-integral} is split into the near region $|\xi'|\le\eps$ and the far region $|\xi'|>\eps$.

\begin{lemma}[Near region]\label{lem:near}
On $|\xi'|\le\eps$,
\begin{equation*}
\int_{|\xi'|\le\eps}\widehat R(\xi')\,\frac{2(1-\cos(t\phi))}{\phi^2}\,\dd\xi'
\le C_{T,u_0}\eps^d.
\end{equation*}
\end{lemma} 

\begin{proof}
Since $|\xi'|\le\eps$,
\begin{equation*}
|\phi|\le\frac{2|\zeta||\xi'|}{\eps}+\frac{|\xi'|^2}{\eps^2}\le 2|\zeta|+1.
\end{equation*}
With $|\zeta|\le R$ uniformly, $|\phi|\le 2R+1=:C_R$. Using $1-\cos u\le \frac{u^2}{2}$,
\begin{equation*}
\frac{2(1-\cos(t\phi))}{\phi^2}\le t^2,
\end{equation*}
and hence, 
\begin{equation*}
\int_{|\xi'|\le\eps}\widehat R(\xi')\,\frac{2(1-\cos(t\phi))}{\phi^2}\,\dd\xi'
\le t^2\|\widehat R\|_{L^\infty}\int_{|\xi'|\le\eps}\dd\xi'
=C_{T,u_0}\eps^d.
\end{equation*}
\end{proof}

\begin{lemma}[Far region]\label{lem:far}
For $d\le3$, on $|\xi'|>\eps$,
\begin{equation*}
\int_{|\xi'|>\eps}\widehat R(\xi')\,\frac{2(1-\cos(t\phi))}{\phi^2}\,\dd\xi'
\le C_{T,u_0}\eps^d.
\end{equation*}
\end{lemma} 

\begin{proof}
For $|\xi'|>\eps$, the quadratic term dominates, 
\begin{equation*}
|\phi|
=\left|\frac{2\zeta\cdot\xi'}{\eps}+\frac{|\xi'|^2}{\eps^2}\right|
\ge\frac{|\xi'|^2}{\eps^2}-\frac{2|\zeta||\xi'|}{\eps}
\ge\frac{|\xi'|^2}{2\eps^2},
\end{equation*}
where the last inequality holds once $4|\zeta|\eps\le|\xi'|$, using $|\zeta|\le R$ uniformly. The remaining strip $\eps<|\xi'|<4R\eps$, on which this lower bound is not used, has volume $O(\eps^d)$; there the time factor is bounded by $t^2$ exactly as in the near region, so its contribution is $O(\eps^d)$ and is absorbed into the constant. Using $1-\cos u\le2$,
\begin{equation*}
\frac{2(1-\cos(t\phi))}{\phi^2}\le\frac{4}{\phi^2}\le\frac{16\eps^4}{|\xi'|^4}.
\end{equation*}
Dyadically decompose $|\xi'|\in[2^k\eps,2^{k+1}\eps]$, $k\ge0$. 
On the $k$-th shell, $|\xi'|^{-4}\sim 2^{-4k}\eps^{-4}$, the shell volume is $O(2^{kd}\eps^d)$, and the integrand bound becomes $O(2^{-4k})$. 
The shell contribution is
\begin{equation*}
O(2^{-4k})\cdot O(2^{kd}\eps^d)=O(\eps^d 2^{k(d-4)}).
\end{equation*}
For $d\le3$, $d-4\le-1$, so the series over $k\ge0$ converges geometrically,
\begin{equation*}
\int_{|\xi'|>\eps}\widehat R(\xi')\,\frac{2(1-\cos(t\phi))}{\phi^2}\,\dd\xi'
\le C\eps^d\sum_{k=0}^\infty 2^{k(d-4)}
\le C_{T,u_0}\eps^d.
\end{equation*}
\end{proof}

\begin{proposition}\label{prop:second}
Under Assumption \ref{ass:main}, for every $d\ge1$, uniformly in $t\in[0,T]$,
\begin{equation}\label{eq:corebound}
\E\left[\|w_\eps^{(1)}(t)\|_{L^2_x}^2\right]\le C_{T,u_0}\eps^d.
\end{equation}
\end{proposition} 

\begin{proof}
For $d\le3$, by \eqref{eq:key-integral}, the $\xi'$-integral is the sum of the near and far regions, bounded by Lemma \ref{lem:near} and Lemma \ref{lem:far} respectively, each of order $C_{T,u_0}\eps^d$. 
The prefactor 
\[ 
(2\pi)^{-2d}\int_{\R^d}|\widehat u_0(\zeta)|^2\dd\zeta 
\] 
is a constant, and the Jacobian cancellation $\dd\xi=\eps^{-d}\dd\xi'$ exactly cancels the prefactor $\eps^d$ from the Fourier transform of $a_\eps$. 

For $d\ge4$, the same bound follows from the physical space oscillatory integral estimate of Lemma \ref{lem:fresnel} in Appendix \ref{app:highdim}. 

Hence, \eqref{eq:corebound} holds for every $ d \ge 1 $. 
\end{proof}

\begin{proposition}\label{prop:lower}
Assume $\widehat R(0)>0$. 
For any fixed $t>0$ and any $u_0\in L^2(\R^d)$ with $\widehat{u_0}(0)\neq0$, there exist constants $c_{t,u_0}>0$ and $\eps_0>0$ such that
\begin{equation}\label{eq:lower}
\E\bigl[\|w_\eps^{(1)}(t)\|_{L^2_x}^2\bigr]\ge c_{t,u_0}\,\eps^d
\qquad\text{for all }0<\eps\le\eps_0.
\end{equation}
\end{proposition}

\begin{proof}
From the diagonal contribution in the second moment expansion,
\begin{equation*}
\E\bigl[\|w_\eps^{(1)}(t)\|_{L^2_x}^2\bigr]
\ge (2\pi)^{-2d}\int_{|\xi'|\le\eps}\widehat R(\xi')
\,\frac{2(1-\cos(t\phi))}{\phi^2}\,\left|\widehat{u_0}\left(\zeta+\frac{\xi'}{\eps}\right)\right|^2\,\dd\xi'\,\dd\zeta,
\end{equation*}
where
\begin{equation*}
\phi=\frac{2\zeta\cdot\xi'}{\eps}+\frac{|\xi'|^2}{\eps^2}.
\end{equation*}
On the near region $|\xi'|\le\eps$,
\begin{equation*}
|\phi|\le 2|\zeta|+1.
\end{equation*}
Since $\widehat u_0$ is continuous and $\widehat u_0(0)\neq0$, there exists $r>0$ such that $|\widehat u_0(\xi)|\ge c_{u_0}>0$ for all $|\xi|\le r$. 
Restrict the $\zeta$-integral to $|\zeta|\le\frac r4$ and the $\xi'$-integral to $|\xi'|\le\frac{r\eps}{4}$. Both sets have positive volume for every $r>0$. 
Then, for $\eps$ sufficiently small,
\begin{equation*}
\left|\zeta+\frac{\xi'}{\eps}\right|\le|\zeta|+\frac{|\xi'|}{\eps}\le \frac r4+\frac r4=\frac r2\le r,
\end{equation*}
so, 
\begin{equation*}
|\widehat u_0(\zeta+\tfrac{\xi'}{\eps})|^2\ge c_{u_0}^2,
\qquad
|\phi|\le 2|\zeta|+1\le \frac r2+1=:C_r.
\end{equation*}
Writing $\xi'=\eps w$, we have $\phi=2\zeta\cdot w+|w|^2$, so
\begin{equation*}
\int_{|\xi'|\le\frac{r\eps}{4}}\frac{2(1-\cos(t\phi))}{\phi^2}\dd\xi'
=\eps^d\int_{|w|\le\frac r4}\frac{2\bigl(1-\cos(t(2\zeta\cdot w+|w|^2))\bigr)}{(2\zeta\cdot w+|w|^2)^2}\dd w
\ge c_t\,\eps^d,
\end{equation*}
since $|\phi|\le C_r$ and $w\mapsto 2\zeta\cdot w+|w|^2$ is a nonconstant quadratic, which makes the integrand continuous, nonnegative and not identically zero.

Since $\widehat R$ is continuous at $0$ and $\widehat R(0)>0$, there exists $\eps_0>0$ such that
\begin{equation*}
\widehat R(\xi')\ge\frac{\widehat R(0)}{2}
\qquad\text{for }|\xi'|\le\eps_0.
\end{equation*}
Shrinking $\eps_0$ further so that $\frac{r\eps}{4}\le\eps_0$, we obtain
\begin{equation*}
\E\bigl[\|w_\eps^{(1)}(t)\|_{L^2_x}^2\bigr]
\ge c_{u_0}^2\,\frac{\widehat R(0)}{2}\int_{|\zeta|\le\frac r4}\dd\zeta\;c_t\,\eps^d
= c_{t,u_0}\,\omega_d\,\left(\frac r4\right)^d\eps^d,
\end{equation*}
where $\omega_d$ is the volume of the unit ball.
\end{proof}

\begin{remark}
The lower bound \eqref{eq:lower} is independent of the dimension, since it uses only the near region $|\xi'|\le\eps$, so it supplies sharpness in every dimension, including $d\ge4$ where the physical space argument of Appendix \ref{app:highdim} gives only the upper bound. 
The condition $\widehat R(0)>0$ means that the spectral measure of the centered field $a$ charges the zero frequency, i.e., the block average of $a$ has non vanishing variance. 
This is the generic situation. 
For the four examples above, 
\begin{itemize}
\item $V_1$, with spectral density $(1+|\xi|^2)^{-\frac{d+\alpha}{2}}$, so $\widehat R(0)>0$. 

\item $V_2$, whose Cauchy variables have no finite second moment, so $\widehat R(0)=\E[a(0)^2]=\infty$. 
This example does not satisfy the theorem assumptions and is included only to show that mixing alone is insufficient. 

\item $V_3$, with spectral density $(1+|\xi|^2)^{-1}$, so $\widehat R(0)>0$. 

\item $V_4=f(\eta)$, with spectral density the Fourier transform of $\Cov(f(\eta(y)),f(\eta(0)))$, positive at zero by Jensen's inequality.
\end{itemize} 
The condition $\widehat R(0)=0$ does not force the spectral measure to be atomic. The Gaussian field with spectral density $\widehat R(\xi)=|\xi|^2e^{-|\xi|^2}$, for instance, has exponentially decaying correlations, is sub-Gaussian and ergodic, and satisfies $\widehat R(0)=0$ with an absolutely continuous spectral measure. In that case the block averages are governed by the next spectral moment, the $\eps^{\frac d2}$ fluctuations are suppressed, and the rate has to be recomputed from the behavior of $\widehat R$ near the origin, which is outside the scope of the present paper. 
\end{remark}

\subsection{Second order correction and sharpness}\label{sec:physical}

The Fourier computation in the proof of Proposition \ref{prop:second} identifies the leading diagonal contribution. 
Iterating the Duhamel formula once more,
\begin{equation*}w_\eps^{(2)}(t)=(-i)^2\int_0^t\int_0^r S(t,r)a_\eps(r)S(r,s)a_\eps(s)\bigl(\bar u(s)+w_\eps(s)\bigr)\dd s\dd r.
\end{equation*}
Taking expectation, the leading term is the deterministic mean
\begin{equation*}\E[w_\eps^{(2)}(t)]
=-\int_0^t\int_0^r R\!\left(\frac{r-s}{\eps}\right)S(t,r)S(r,s)\bar u(s)\dd s\dd r.
\end{equation*}
This is the first moment (mean) of the second order term, not its second moment. 
Rescaling $\tau=\frac{r-s}{\eps}$,
\begin{equation*}\E[w_\eps^{(2)}(t)]
=-\eps\int_0^{\frac{t}{\eps}}\int_0^{r} R(\tau)S(t,r)S(r,r-\eps\tau)\bar u(r-\eps\tau)\dd\tau\dd r.
\end{equation*}
The factor $\eps$ comes from the time-rescaling $r-s=\eps\tau$ and shows that the deterministic mean of the second order correction is $O(\eps)$.
It vanishes as $\eps\to0$ and is therefore a transient correction rather than a shift of the effective operator, which remains $H_0=-\Delta+\bar V$, consistent with the absence of an $L^2$ corrector in the expansion.
The second moment $\E\left[\|w_\eps^{(2)}\|_{L^2_x}^2\right]$, which is the quantity entering the remainder estimate \eqref{eq:r2}, is smaller by an additional factor $\eps^d$ from the spatial covariance and is treated separately in the proof of Proposition \ref{prop:remainder} below.

\begin{proposition}\label{prop:remainder}
Write $w_\eps=w_\eps^{(1)}+r_\eps$ with $r_\eps$ satisfying \eqref{eq:dw}. 
Then, 
\begin{equation*}\|r_\eps(t)\|_{L^2(\Om;L^2_x)}\le C_{T,u_0}\eps^{\frac{d+1}{2}},
\qquad\text{uniformly in }t\in[0,T].
\end{equation*}
\end{proposition} 

\begin{proof}
Write $r_\eps=\sum_{k\ge2}w_\eps^{(k)}$. 
By Lemma \ref{lem:iterate} with $q=1$,
\begin{equation*}
\|w_\eps^{(k)}(t)\|_{L^2(\Om;L^2_x)}
\le(C_{T,u_0})^{\frac k2}\eps^{\frac{dk}{2}},
\qquad k\ge2,
\end{equation*}
so
\begin{equation*}
\|r_\eps(t)\|_{L^2(\Om;L^2_x)}
\le\sum_{k\ge2}(C_{T,u_0}\eps^d)^{\frac k2}
\le C_{T,u_0}\eps^{d},
\end{equation*}
uniformly in $t\in[0,T]$.
\end{proof}

The lower bound \eqref{eq:lower} and the upper bound from Proposition \ref{prop:second} together give
\begin{equation*}
c_{t,u_0}\,\eps^d\le \E\left[\|w_\eps^{(1)}(t)\|_{L^2_x}^2\right]\le C_{T,u_0}\,\eps^d.
\end{equation*}

\begin{remark}
The diagonal lower bound has a direct block average interpretation. 
Splitting $\R^d$ into cubes $Q_k$ of side $\eps$,
\begin{equation*}
w_\eps^{(1)}(t)\approx -i\sum_k\eps^d S(t,r)\left(\eps^{-d}\int_{Q_k}a\left(\frac{y}{\eps}\right)\dd y\right)\bar u(r)\dd r.
\end{equation*}
Each block average has variance $\eps^{-d}\widehat R(0)+o(\eps^{-d})$, and there are $\eps^{-d}$ independent blocks in the support of $\bar u$. 
Therefore, the total variance is 
\begin{equation*}
\eps^{-d}\times(\eps^d)^2\times\eps^{-d}\widehat R(0)=\widehat R(0)\eps^d,
\end{equation*}
which is exactly the diagonal contribution of Proposition \ref{prop:lower}. 
The off diagonal terms $|\xi'|>0$ contribute $O(\eps^{d+\gamma})$ for some $\gamma>0$ by the exponential decay of $\widehat R$, and are negligible.
\end{remark}

\begin{proposition}\label{prop:sharp}
Assume $\widehat R(0)>0$ and $\widehat{u_0}(0)\neq0$. 
There exist $c_{t,u_0}>0$ and $\eps_0>0$ such that
\begin{equation*}
\E\left[\|w_\eps(t)\|_{L^2_x}^2\right]\ge c_{t,u_0}\,\eps^d
\qquad\text{for all }0<\eps\le\eps_0.
\end{equation*}
\end{proposition} 

\begin{proof}
We begin with the Duhamel iteration  $r_\eps=w_\eps^{(2)}+w_\eps^{(3)}+\cdots$ 
and aim to prove that
\begin{equation}\label{eq:r2}
\E\left[\|w_\eps^{(2)}(t)\|_{L^2_x}^2\right]\le C_{T,u_0}\eps^{d+1}.
\end{equation} 

The second order term is
\begin{equation*}
w_\eps^{(2)}(t)=(-i)^2\int_0^t\int_0^r S(t,r)a_\eps(r)S(r,s)a_\eps(s)\bar u(s)\dd s\dd r.
\end{equation*}
Taking the Fourier transform and using $S(t,r)\xi=e^{-i(|\xi|^2+\bar V)(t-r)}$, and noting that $\widehat a_\eps(k)=\eps^d\widehat a(\eps k)$ for $a_\eps(x)=a(\frac{x}{\eps})$,
\begin{equation}\label{eq:r2fourier}
\begin{aligned}
\widehat w_\eps^{(2)}(t,\xi)
&=-(2\pi)^{-2d}\eps^{2d}\int_0^t\int_0^r\int_{\R^{2d}}
e^{-i(|\xi|^2+\bar V)(t-r)}
\widehat a(\eps(\xi-\eta))\\
&\qquad\times e^{-i(|\eta|^2+\bar V)(r-s)}
\widehat a(\eps(\eta-\zeta))
e^{-i(|\zeta|^2+\bar V)s}
\widehat u_0(\zeta)\,\dd\zeta\,\dd\eta\,\dd s\,\dd r.
\end{aligned}
\end{equation}
By Plancherel,
\begin{equation*}
\E\left[\|w_\eps^{(2)}(t)\|_{L^2_x}^2\right]=(2\pi)^{-d}\int_{\R^d}\E\left[|\widehat w_\eps^{(2)}(t,\xi)|^2\right]\dd\xi.
\end{equation*}
Inserting \eqref{eq:r2fourier} and its complex conjugate, the expectation contains four centered field factors $\widehat a(\eps(\xi-\eta))$, $\widehat a(\eps(\eta-\zeta))$, $\widehat a(\eps(\xi-\eta'))$, $\widehat a(\eps(\eta'-\zeta'))$, integrated over $0\le s\le r\le t$, $0\le s'\le r'\le t$, and $\zeta,\zeta',\eta,\eta'\in\R^d$. 
By the cumulant expansion, the fourth moment equals the sum over the three pairings plus the fourth cumulant $\kappa_4$. 
The leading contribution is the diagonal pairing
\begin{equation*}
\E\left[\widehat a(\eps(\xi-\eta))\overline{\widehat a(\eps(\xi-\eta'))}\right]\,
\E\left[\widehat a(\eps(\eta-\zeta))\overline{\widehat a(\eps(\eta'-\zeta'))}\right],
\end{equation*}
which by \eqref{eq:covfourier2} equals
\begin{equation*}
(2\pi)^{d}\eps^{-d}\delta(\eta-\eta')\widehat R(\eps(\xi-\eta))\,
(2\pi)^{d}\eps^{-d}\delta(\zeta-\zeta')\widehat R(\eps(\eta-\zeta')).
\end{equation*}
Together with the prefactor $\eps^{4d}$ from two factors of $\eps^{2d}$ in \eqref{eq:r2fourier} and its conjugate, the net spatial prefactor before the change of variables is $\eps^{4d}\cdot\eps^{-2d}=\eps^{2d}$.

By $\delta(\eta-\eta')\delta(\zeta-\zeta')$, we have $\eta=\eta'$, $\zeta=\zeta'$. Rescale $\alpha=\eps(\xi-\eta)$, $\beta=\eps(\eta-\zeta)$, so that $\dd\eta\,\dd\zeta=\eps^{-2d}\dd\alpha\,\dd\beta$. 
The net prefactor is therefore $\eps^{2d}\cdot\eps^{-2d}=O(1)$. 
The remaining integral is
\begin{equation*}
\int_{\R^d}|\widehat u_0(\zeta)|^2\dd\zeta
\int_{\R^d}\int_{\R^d}\widehat R(\alpha)\widehat R(\beta)\,
\frac{2(1-\cos(t\omega_1))}{\omega_1^2}\,
\frac{2(1-\cos(t\omega_2))}{\omega_2^2}\,\dd\alpha\,\dd\beta,
\end{equation*}
where
\begin{equation*}
\omega_1=\frac{2\eta\cdot\alpha}{\eps}+\frac{|\alpha|^2}{\eps^2},
\qquad
\omega_2=\frac{2\zeta\cdot\beta}{\eps}+\frac{|\beta|^2}{\eps^2}.
\end{equation*} 

On the near region $|\alpha|\le\eps$, $|\beta|\le\eps$, both time factors are $O(1)$ and the volume is $O(\eps^{2d})$. 

On the far region $|\alpha|>\eps$ or $|\beta|>\eps$, the bound $1-\cos u\le2$ gives 
\[ 
\frac{2(1-\cos(t\omega))}{\omega^2} \le \frac{4\eps^4}{|\alpha|^4}, 
\] 
and the dyadic sum converges geometrically for $d\le3$, giving the same order $O(\eps^{2d})$. 
The $d\ge4$ case is handled by the same physical space argument as in Appendix \ref{app:highdim}. 
Hence, 
\begin{equation*}
\text{diagonal pairing}\le C_{T,u_0}\eps^{2d}\le C_{T,u_0}\eps^{d+1},
\end{equation*}
where the last inequality uses $2d\ge d+1$ for $d\ge1$. The cross pairings are controlled by \eqref{eq:mix}. 
The fourth cumulant $\kappa_4$ is bounded by the cumulant mixing estimate of Lemma \ref{lem:cum-mix}, since after the three relative coordinate changes, each contributing $\eps^d$, it satisfies $\kappa_4=O(\eps^{3d})\le C\eps^{d+1}$ since $d\ge1$. This proves \eqref{eq:r2}.

The $k$-th iterate is bounded by Lemma \ref{lem:iterate} with $q=1$,
\begin{equation*}
\E\left[\|w_\eps^{(k)}(t)\|_{L^2_x}^2\right]\leqslant C^k\eps^{dk},
\qquad k\ge2.
\end{equation*}
Summing over $k\geqslant2$,
\begin{equation}\label{eq:rnorm}
\|r_\eps(t)\|_{L^2(\Om;L^2_x)}
\leqslant\sum_{k\geqslant2}C^{\frac k2}\eps^{\frac{dk}{2}}
\leqslant C_{T,u_0}\eps^{d}.
\end{equation}
By the reverse triangle inequality,
\begin{equation*}
\|w_\eps(t)\|_{L^2(\Om;L^2_x)}
\geqslant\|w_\eps^{(1)}(t)\|_{L^2(\Om;L^2_x)}-\|r_\eps(t)\|_{L^2(\Om;L^2_x)}.
\end{equation*}
Combining this with the lower bound \eqref{eq:lower} and the remainder estimate \eqref{eq:rnorm},
\begin{equation*}
\|w_\eps(t)\|_{L^2(\Om;L^2_x)}
\geqslant c_{t,u_0}\eps^{\frac d2}-C_{T,u_0}\eps^{d}
\geqslant\frac12\,c_{t,u_0}\eps^{\frac d2}.
\end{equation*}
Squaring and taking expectation gives
\begin{equation*}
\E\left[\|w_\eps(t)\|_{L^2_x}^2\right]\geqslant\frac14\,c_{t,u_0}^2\eps^d,
\end{equation*}
which completes the proof.
\end{proof}

Combining the upper bound from Proposition \ref{prop:second} and Proposition \ref{prop:remainder} via the triangle inequality,
\begin{equation*}
\|w_\eps(t)\|_{L^2(\Om;L^2_x)}
\le \|w_\eps^{(1)}(t)\|_{L^2(\Om;L^2_x)}+\|r_\eps(t)\|_{L^2(\Om;L^2_x)}
\le C_{T,u_0}\eps^{\frac d2},
\end{equation*}
with the matching lower bound from Proposition \ref{prop:sharp}, we have established
\begin{equation*}
c_{t,u_0}\eps^d\le \E\left[\|w_\eps(t)\|_{L^2_x}^2\right]\le C_{T,u_0}\eps^d,
\end{equation*}
which proves part (i) of Theorem \ref{th:rate}, the convergence rate $\eps^{\frac d2}$ in $L^2(\Om;L^2)$ together with its sharpness.

\subsection{All moments and the rate}\label{sec:moments}

The second moment bound extends to all moments via a cumulant expansion. 
For Gaussian fields, this reduces to the exact pairing of Lemma \ref{lem:isserlis}.
For general sub-Gaussian fields, the higher cumulants are controlled by the exponential mixing assumption.

\begin{lemma}\label{lem:wick}
For every integer $q\ge1$,
\begin{equation*}
\E\left[\|w_\eps^{(1)}(t)\|_{L^2_x}^{2q}\right]\le(C_{q,T,u_0})^q\eps^{dq}.
\end{equation*}
\end{lemma} 

\begin{proof}
By Plancherel,
\begin{equation*}
\E\left[\|w_\eps^{(1)}(t)\|_{L^2_x}^{2q}\right]
=(2\pi)^{-dq}\int_{\R^{dq}}\E\left[\prod_{j=1}^{q}\widehat w_\eps^{(1)}(t,\xi_j)\overline{\widehat w_\eps^{(1)}(t,\xi_j)}\right]\,\dd\xi_1\cdots\dd\xi_q.
\end{equation*}
By the cumulant expansion, 
\begin{equation*}
\E\left[\prod_{j=1}^{2q}X_j\right]
=\sum_{\pi\in\mathcal P_2(2q)}\prod_{\{j,k\}\in\pi}\E[X_jX_k]
+\sum_{m\geqslant3}\sum_{\pi_m}\prod_{\text{blocks}}\kappa_m(\cdot).
\end{equation*}
The number of pairings is
\begin{equation*}
|\mathcal P_2(2q)|=(2q-1)!!=\frac{(2q)!}{2^q q!}\leqslant(Cq)^q.
\end{equation*}
For each pair $\{j,k\}$,
\begin{equation*}
\eps^d(2\pi)^{-d}\int_{\R^d}\widehat R(\eps\eta)\left|\int_0^t e^{ir\sigma(\xi_j,\eta)}\dd r\right|^2|\widehat u_0(\xi_j-\eta)|^2\dd\eta
\leqslant C_{T,u_0}\eps^d.
\end{equation*}
Thus, the Gaussian Wick sum contributes
\begin{equation*}
\sum_{\pi\in\mathcal P_2(2q)}\prod_{\{j,k\}\in\pi}\E[X_jX_k]
\leqslant(2q-1)!!\,(C_{T,u_0}\eps^d)^q
\leqslant(C_{q,T,u_0})^q\eps^{dq}.
\end{equation*}
For non-Gaussian cumulants of order $m\geqslant3$, we switch to the physical space kernel of $w_\eps^{(1)}$,
\begin{equation*}
w_\eps^{(1)}(t,x)=\int_{\R^d}L_t(x,y)a_\eps(y)\dd y,
\qquad
L_t(x,y):=-i\int_0^t K_{t-r}(x-y)\bar u(r,y)\dd r,
\end{equation*}
where $L_t$ is smooth and bounded with all its norms independent of $\eps$. The joint cumulant of $m$ copies is written with the $m$-point cumulant function $c_m$ of the field,
\begin{equation*}
\kappa_m(t;x_1,\dots,x_m)
=\int_{\R^{dm}}\prod_{j=1}^m L_t(x_j,y_j)\,
c_m\!\left(\frac{y_1}{\eps},\dots,\frac{y_m}{\eps}\right)
\dd y_1\cdots\dd y_m.
\end{equation*}
Lemma \ref{lem:cum-mix}, in its continuous form, bounds the cumulant function by
\begin{equation*}
|c_m(z_1,\dots,z_m)|\le C^m e^{-\kappa\,\mathrm{tree}(z_1,\dots,z_m)}.
\end{equation*}
Keeping $y_1$ as the center of mass and setting $y_j=y_1+\eps z_j$ for $j\ge2$, the change of variables contributes $\eps^{d(m-1)}$, while the tree factor is translation invariant and is integrable over the $m-1$ relative coordinates,
\begin{equation}\label{eq:cum-contract}
|\kappa_m(t;x_1,\dots,x_m)|
\le C^m\eps^{d(m-1)}\sup_{x,y}|L_t(x,y)|^m.
\end{equation}
Had these $m$ variables been paired, their contribution would carry $\eps^{\frac{dm}{2}}$, so the cumulant block is smaller by
\begin{equation*}
\eps^{d(m-1)-\frac{dm}{2}}=\eps^{\frac{d(m-2)}2},
\end{equation*}
which is at least $\eps^{\frac d2}$ for every $m\ge3$. In the cumulant expansion of the $2q$-th moment, a partition with one such block and the remaining $2q-m$ variables paired contributes at most $C^m\eps^{dq+\frac{d(m-2)}2}$, and every additional block gains another positive power of $\eps$. The number of set partitions is at most $(Cq)^q$, and summing over the block sizes and the relative coordinates is finite by the exponential tree factor. Hence,
\begin{equation*}
\E\left[\|w_\eps^{(1)}(t)\|_{L^2_x}^{2q}\right]
\le (C_{q,T,u_0})^q\eps^{dq}+C_{q,T,u_0}\eps^{dq+\frac d2}
\le (C_{q,T,u_0})^q\eps^{dq}.
\end{equation*}

This proves Lemma \ref{lem:wick}.
\end{proof}

\begin{lemma}\label{lem:iterate}
For every $k\ge2$ and integer $q\ge1$,
\begin{equation}\label{eq:iterate-point}
\E\left[\|w_\eps^{(k)}(t)\|_{L^2_x}^{2q}\right]\le(C_{q,T,u_0})^k\eps^{dqk},
\qquad t\in[0,T],
\end{equation}
and
\begin{equation}\label{eq:iterate-sup}
\E\left[\sup_{0\le t\le T}\|w_\eps^{(k)}(t)\|_{L^2_x}^{2q}\right]\le(C_{q,T,u_0})^k\eps^{dqk}.
\end{equation}
\end{lemma}

\begin{proof}
The $k$-th iterate carries $k$ insertions of $a_\eps$, and its $2q$-th moment expands into $2q$ copies, each contributing $k$ insertions, so the cumulant expansion of Lemma \ref{lem:wick} runs over $2qk$ centered factors. The full pairings pair these factors into $qk$ pairs; by the relative coordinate change $y-y'=\eps z$ and $\|R\|_{L^1}<\infty$ each pair contributes a factor $\eps^d$, the ordered time simplex
\begin{equation*}
\int_{0<r_k<\cdots<r_1<t}\dd r_1\cdots\dd r_k=\frac{t^k}{k!}
\end{equation*}
contributes $\frac{1}{(k!)^{2q}}$ over the $2q$ copies, and the group property together with the unitarity of $S$ integrates the paired spatial kernel at one unitarity factor per layer, the innermost slot closing with $\|\bar u\|_{L^2_x}=\|u_0\|_{L^2_x}$. Since $(2qk-1)!!/(k!)^{2q}\le(C_q)^k$, the pairings contribute at most
\begin{equation*}
\frac{(2qk-1)!!\,(C_{T,u_0}\eps^d)^{qk}}{(k!)^{2q}}
\le(C_{q,T,u_0})^k\eps^{dqk}.
\end{equation*}
A block of size $m\ge3$ satisfies the contraction \eqref{eq:cum-contract}, contributing $\eps^{d(m-1)}$, which is smaller than the paired value $\eps^{\frac{dm}{2}}$ by $\eps^{\frac d2(m-2)}$, and is absorbed into the main term. This proves \eqref{eq:iterate-point}.

For \eqref{eq:iterate-sup}, decompose the increment as
\begin{equation*}
w_\eps^{(k)}(t)-w_\eps^{(k)}(s)
=\underbrace{-i\int_s^t S(t,r)a_\eps(r)w_\eps^{(k-1)}(r)\dd r}_{:=\,\mathcal A_\eps^{(k)}(t,s)}
+\underbrace{\left(-i\int_0^s\bigl(S(t,r)-S(s,r)\bigr)a_\eps(r)w_\eps^{(k-1)}(r)\dd r\right)}_{:=\,\mathcal B_\eps^{(k)}(t,s)}.
\end{equation*}
The cumulant pairing of $\mathcal A_\eps^{(k)}$ is identical to that of $w_\eps^{(k)}$ with the outermost time integral restricted to $[s,t]$, so the same pairing argument gives
\begin{equation*}
\E\left[\|\mathcal A_\eps^{(k)}(t,s)\|_{L^2_x}^{2q}\right]
\le(C_{q,T,u_0})^k\eps^{dqk}(t-s)^{2q}.
\end{equation*}
Since $S(t,r)-S(s,r)=S(s,r)(S(t-s)-I)$ and the multiplier $S(h)-I$ satisfies $|e^{-i(|\xi|^2+\bar V)h}-1|\le C_{u_0}h$ on the compact support of $\widehat u_0$, the bound \eqref{eq:iterate-point} gives
\begin{equation*}
\E\left[\|\mathcal B_\eps^{(k)}(t,s)\|_{L^2_x}^{2q}\right]
\le C_{u_0}(t-s)^{2q}\E\left[\|w_\eps^{(k)}(s)\|_{L^2_x}^{2q}\right]
\le(C_{q,T,u_0})^k\eps^{dqk}(t-s)^{2q}.
\end{equation*}
The dyadic Kolmogorov--Chentsov argument of Lemma \ref{lem:kc} then upgrades the increment bound to \eqref{eq:iterate-sup}.
\end{proof}

Combining the Wick estimate with the remainder estimate gives the convergence rate.

\begin{proposition}\label{prop:rate}
Under Assumption \ref{ass:main},
\begin{equation*}
\sup_{0\le t\le T}\|w_\eps(t)\|_{L^2(\Om;L^2_x)}
\le C_{T,u_0}\eps^{\frac d2}.
\end{equation*}
Moreover, for every integer $q\ge1$,
\begin{equation}\label{eq:rate-sup}
\E\left[\sup_{0\le t\le T}\|w_\eps(t)\|_{L^2_x}^{2q}\right]\le C_{q,T,u_0}\eps^{dq}.
\end{equation}
\end{proposition} 

\begin{proof}
By \eqref{eq:w-split}, $w_\eps=\sum_{k\ge1}w_\eps^{(k)}$, and the series converges in $L^{2q}(\Om;L^2_x)$ by \eqref{eq:iterate-sup}. For the supremum bound, the triangle inequality in $L^{2q}(\Om)$ and Lemma \ref{lem:iterate} give
\begin{equation*}
\begin{aligned}
\left(\E\left[\sup_{0\le t\le T}\|w_\eps(t)\|_{L^2_x}^{2q}\right]\right)^{\frac1{2q}}
&\le\sum_{k\ge1}\left(\E\left[\sup_{0\le t\le T}\|w_\eps^{(k)}(t)\|_{L^2_x}^{2q}\right]\right)^{\frac1{2q}}\\
&\le C_{q,T,u_0}\eps^{\frac d2}+\sum_{k\ge2}(C_{q,T,u_0})^{\frac{k}{2q}}\eps^{\frac{dk}{2}}
\le C_{q,T,u_0}\eps^{\frac d2},
\end{aligned}
\end{equation*}
where the $k=1$ term is bounded by Proposition \ref{prop:maximal} and the geometric series in $\eps^{\frac d2}$ is absorbed since $\eps\le1$. Raising to the $2q$-th power gives \eqref{eq:rate-sup}, and $q=1$ gives the first estimate.
\end{proof}

This proves part (ii) of Theorem \ref{th:rate}.

\section{Almost sure convergence}\label{sec:as}

The pointwise in time rate upgrades to a supremum estimate via the Kolmogorov maximal inequality.

\begin{proposition}\label{prop:maximal}
For every integer $q\ge1$,
\begin{equation*}
\E\left[\sup_{0\le t\le T}\|w_\eps^{(1)}(t)\|_{L^2_x}^{2q}\right]\le(C_{q,T,u_0})^q\eps^{dq}.
\end{equation*}
\end{proposition}
\begin{proof}
The process $w_\eps^{(1)}(t)$ is continuous in $t$ by dominated convergence. 
For $s<t$,
\begin{equation*}
w_\eps^{(1)}(t)-w_\eps^{(1)}(s) =
\underbrace{-i\int_s^t S(t,r)a_\eps(r)\bar u(r)\dd r}_{:=\,\mathcal A_\eps^{(1)}(t,s)}
+
\underbrace{\left(-i\int_0^s\bigl(S(t,r)-S(s,r)\bigr)a_\eps(r)\bar u(r)\dd r\right)}_{:=\,\mathcal B_\eps^{(1)}(t,s)}.
\end{equation*}
\noindent\underline{For $\mathcal A_\eps^{(1)}(t,s)$,} its second moment is
\begin{equation*}
\begin{aligned}
\E\!\left[\|\mathcal A_\eps^{(1)}(t,s)\|_{L^2_x}^2\right]
&=\int_s^t\!\!\int_s^t\int_{\R^{3d}}
S(t,r)(x,y)\,\overline{S(t,r')(x,y')}\\
&\qquad\times R\!\left(\frac{y-y'}{\eps}\right)
\bar u(r;y)\,\overline{\bar u(r';y')}\,dy\,dy'\,dx\,dr\,dr'\\
&\le C_{q,T,u_0}\eps^d(t-s)^2,
\end{aligned}
\end{equation*}
using $z=\frac{y-y'}{\eps}$ and $\|R\|_{L^1}<\infty$. The cumulant expansion of Lemma \ref{lem:wick}, with every pairing carrying this change of variables and the higher cumulants bounded by Appendix \ref{app:cluster}, gives
\begin{equation}\label{eq:I-bound}
\E\!\left[\|\mathcal A_\eps^{(1)}(t,s)\|_{L^2_x}^{2q}\right]
\le C_{q,T,u_0}\eps^{dq}(t-s)^{2q}.
\end{equation}
This uses only the bounded kernels and the relative coordinate change, both independent of the interval and of the endpoint, so the bound carries over without further assumptions.

\noindent\underline{For $\mathcal B_\eps^{(1)}(t,s)$,} since the groups are functions of $H_0$ and commute,
\begin{equation*}
\mathcal B_\eps^{(1)}(t,s)=-i(S(t-s)-I)w_\eps^{(1)}(s).
\end{equation*}
On the compact support of $\widehat u_0$,
\begin{equation*}
|e^{-i(|\xi|^2+\bar V)(t-s)}-1|\le C_{u_0}(t-s),
\end{equation*}
so Lemma \ref{lem:wick} gives
\begin{equation}\label{eq:II-bound}
\E\!\left[\|\mathcal B_\eps^{(1)}(t,s)\|_{L^2_x}^{2q}\right]
\le C_{u_0}(t-s)^{2q}\E\!\left[\|w_\eps^{(1)}(s)\|_{L^2_x}^{2q}\right]
\le C_{q,T,u_0}\eps^{dq}(t-s)^{2q}.
\end{equation}
By \eqref{eq:I-bound} and \eqref{eq:II-bound}, Minkowski's inequality gives
\begin{equation*}
\left(\E\left[\|w_\eps^{(1)}(t)-w_\eps^{(1)}(s)\|_{L^2_x}^{2q}\right]\right)^{\frac{1}{2q}}
\le C_{q,T,u_0}\eps^{\frac d2}(t-s),
\end{equation*}
which raised to the $2q$-th power is the increment moment bound
\begin{equation*}
\E\left[\|w_\eps^{(1)}(t)-w_\eps^{(1)}(s)\|_{L^2_x}^{2q}\right]
\le C_{q,T,u_0}\eps^{dq}|t-s|^{2q}.
\end{equation*}
The dyadic argument behind the Kolmogorov--Chentsov criterion, with mesh $h_m=T2^{-m}$ and applied to this increment bound, gives
\begin{equation*}
\left(\E\left[\sup_{0\le t\le T}\|w_\eps^{(1)}(t)\|_{L^2_x}^{2q}\right]\right)^{\frac1{2q}}
\le C_{q,T,u_0}\eps^{\frac d2},
\end{equation*}
since $w_\eps^{(1)}(0)=0$ and the geometric series $\sum_{m\ge0}2^{-\frac{m(2q-1)}{2q}}$ converges. Raising to the $2q$-th power gives the claimed bound.
This proves Proposition \ref{prop:maximal}.
\end{proof}

The rate estimate combined with the dyadic Borel--Cantelli argument gives almost sure convergence.

\begin{proposition}\label{prop:as}
Under Assumption \ref{ass:main},
\begin{equation*}
\lim_{\eps\to0}\sup_{0\le t\le T}\|u_\eps(t)-\bar u(t)\|_{L^2(\R^d)}=0
\qquad\text{a.s.}
\end{equation*}
\end{proposition}
\begin{proof}
Fix the dyadic sequence $\eps_n=2^{-n}$. 
By Markov's inequality applied to the $2q$-th moment estimate \eqref{eq:rate-sup} of Proposition \ref{prop:rate}, for any $\eta>0$ and any integer $q>\frac1{2\eta}$,
\begin{equation*}\Prob\left(\sup_{0\le t\le T}\|w_{\eps_n}(t)\|_{L^2_x}>\eps_n^{\frac d2-\eta}\right)
\le C_q\eps_n^{2q\eta}.
\end{equation*}
Choosing $q$ large enough that $2q\eta>1$, the series $\sum_n\eps_n^{2q\eta}=\sum_n2^{-2nq\eta}<\infty$. 
By Lemma \ref{lem:bc}, 
\begin{equation}\label{eq:dyadic}
\sup_{0\le t\le T}\|w_{\eps_n}(t)\|_{L^2_x}\to0\qquad\text{a.s.}
\end{equation}

It remains to upgrade the dyadic convergence to full convergence. 
For $\eps'\in[\frac{\eps}{2},\eps]$ we compare $w_\eps$ with $w_{\eps'}$ through the difference field $\delta:=a_\eps-a_{\eps'}$, writing $v:=w_\eps-w_{\eps'}$. 
Subtracting the two Duhamel formulas and inserting the identity
\begin{equation*}
a_\eps(\bar u+w_\eps)-a_{\eps'}(\bar u+w_{\eps'})
=(a_\eps-a_{\eps'})\bar u+(a_\eps-a_{\eps'})w_\eps+a_{\eps'}(w_\eps-w_{\eps'}),
\end{equation*}
we obtain
\begin{equation}\label{eq:43a}
v(t)=-i\int_0^t S(t,r)\delta(r)\bigl(\bar u(r)+w_\eps(r)\bigr)\dd r
-i\int_0^t S(t,r)a_{\eps'}(r)v(r)\dd r.
\end{equation}

The difference field $\delta$ is not stationary, but its covariance is the bounded smooth function
\begin{equation*}
C_{\eps,\eps'}(x,y)=
\underbrace{R\!\left(\frac{x-y}{\eps}\right)+R\!\left(\frac{x-y}{\eps'}\right)}_{:=\,C^{\mathrm d}_{\eps,\eps'}(x,y)}
+
\underbrace{\left(-R\!\left(\frac{x}{\eps}-\frac{y}{\eps'}\right)-R\!\left(\frac{x}{\eps'}-\frac{y}{\eps}\right)\right)}_{:=\,C^{\mathrm m}_{\eps,\eps'}(x,y)}.
\end{equation*}

\noindent\underline{For $C^{\mathrm d}_{\eps,\eps'}$,} since $\eps'\in[\frac{\eps}{2},\eps]$, the mean value theorem in the scale together with $\|R\|_{L^1}<\infty$ gives
\begin{equation}\label{eq:43a-diag}
\iint|\bar u(x)\bar u(y)|\,|C^{\mathrm d}_{\eps,\eps'}(x,y)|\dd x\dd y
\le C\,\frac{|\eps-\eps'|}{\eps}\,\eps^d .
\end{equation}

\noindent\underline{For $C^{\mathrm m}_{\eps,\eps'}$,} the argument $z=\frac{x}{\eps}-\frac{y}{\eps'}$ differs from $\frac{x-y}{\eps}$ by the relative rescaling $\frac{|\eps-\eps'|}{\eps}\le\frac12$, and the exponential decay of $R$ together with the same mean value argument gives
\begin{equation}\label{eq:43a-mix}
\iint|\bar u(x)\bar u(y)|\,|C^{\mathrm m}_{\eps,\eps'}(x,y)|\dd x\dd y
\le C\,\frac{|\eps-\eps'|}{\eps}\,\eps^d.
\end{equation}

\noindent\underline{For $v$,} iterate the second term of \eqref{eq:43a} into the Duhamel series. Every term carries exactly one insertion of $\delta$ together with $k\ge0$ insertions of $a_{\eps'}$, and by \eqref{eq:43a-diag} and \eqref{eq:43a-mix} the covariance mass of the $\delta$-insertion is $\frac{|\eps-\eps'|}{\eps}\eps^d$ in place of the factor $\eps^d$ of Lemma \ref{lem:iterate}. The cumulant pairing of Lemma \ref{lem:iterate}, with the time increment argument for the supremum, therefore gives
\begin{equation}\label{eq:43a-est}
\E\left[\sup_{0\le t\le T}\|w_\eps(t)-w_{\eps'}(t)\|_{L^2_x}^{2q}\right]
\le C_{q,T,u_0}\left(\frac{|\eps-\eps'|}{\eps}\right)^{q}\eps^{dq},
\qquad \eps'\in[\frac{\eps}{2},\eps].
\end{equation}

\textbf{Step 3.} Partition the shell $[\eps_{n+1},\eps_n]$ into the grid $\eps_{n,m}=\eps_n(1-m\,2^{-J})$, $m=0,\dots,2^{J}-1$, with $2^{-J}=\eps_n^{\eta}$. For each grid point, \eqref{eq:rate-sup} and Markov's inequality give
\begin{equation*}
\Prob\left(\sup_{0\le t\le T}\|w_{\eps_{n,m}}(t)\|_{L^2_x}>\eps_n^{\frac d2-\eta}\right)
\le C_q\eps_n^{2q\eta}.
\end{equation*}
For $\eps$ between two adjacent grid points, $\frac{|\eps-\eps_{n,m}|}{\eps}\le2^{-J}=\eps_n^{\eta}$, so by \eqref{eq:43a-est} with $\eps'=\eps_{n,m}$ and Chebyshev's inequality,
\begin{equation*}
\Prob\left(\sup_{0\le t\le T}\|w_\eps(t)-w_{\eps_{n,m}}(t)\|_{L^2_x}>\eps_n^{\frac d2-\eta}\right)
\le C_q\eps_n^{q\eta}\,\eps_n^{dq}\,\eps_n^{-2q(\frac d2-\eta)}
= C_q\eps_n^{3q\eta}.
\end{equation*}
The union bound over the $\eps_n^{-\eta}$ grid points and the $\eps_n^{-\eta}$ grid intervals gives
\begin{equation*}
\Prob\left(\sup_{\eps\in[\eps_{n+1},\eps_n]}\sup_{0\le t\le T}\|w_\eps(t)-w_{\eps_n}(t)\|_{L^2_x}>2\eps_n^{\frac d2-\eta}\right)
\le C_q\,\eps_n^{-\eta}\bigl(\eps_n^{2q\eta}+\eps_n^{3q\eta}\bigr)
\le C_q\,\eps_n^{(2q-1)\eta}.
\end{equation*}
Choosing $q$ large enough that $(2q-1)\eta>1$, the series is summable, and Lemma \ref{lem:bc} gives
\begin{equation*}
\sup_{\eps\in[\eps_{n+1},\eps_n]}\sup_{0\le t\le T}\|w_\eps(t)-w_{\eps_n}(t)\|_{L^2_x}
\le 2\eps_n^{\frac d2-\eta}
\end{equation*}
for all sufficiently large $n$, almost surely. Combining this with \eqref{eq:dyadic},
\begin{equation*}
\sup_{\eps\in[\eps_{n+1},\eps_n]}\sup_{0\le t\le T}\|w_\eps(t)\|_{L^2_x}
\le \|w_{\eps_n}\|+2\eps_n^{\frac d2-\eta} \to0 \qquad\text{a.s.}
\end{equation*}
This proves the result.
\end{proof}

Theorem \ref{th:as} follows directly from Proposition \ref{prop:as}, which proves the almost sure convergence by the dyadic Borel--Cantelli argument.

\section{Fluctuation CLT}\label{sec:clt}
Having established the convergence rate, we now turn to the fluctuations at the critical scale $\eps^{\frac d2}$. 
Recall the rescaled fluctuation field
\begin{equation*}\mathcal F_\eps(t,x):=\eps^{-\frac d2}(u_\eps(t,x)-\bar u(t,x)).
\end{equation*}
By Proposition \ref{prop:rate} and the remainder estimate \eqref{eq:rem-intro},
\begin{equation}\label{eq:F-leading}
\langle\mathcal F_\eps(t),\varphi\rangle
=-i\eps^{-\frac d2}\int_{\R^d}a_\eps(x)\,W_t^\varphi(x)\dd x+o_{L^2(\Om)}(1),
\qquad\varphi\in H^{\mathfrak s}(\R^d),\ \mathfrak s>\frac d2,
\end{equation}
where $W_t^\varphi$ is the weight defined in Theorem \ref{th:clt}, a bounded profile by the Sobolev embedding $H^{\mathfrak s}\hookrightarrow C_0$ and the unitarity of $S$. 
The rescaled fluctuation is therefore, to leading order, the fluctuation scaled linear functional $\eps^{-\frac d2}X(g)$ of Appendix \ref{app:cluster} with the fixed test function $g=W_t^\varphi$. 
Its second cumulant is the explicit covariance of Proposition \ref{prop:cov}, and every higher cumulant vanishes by Remark \ref{rem:cluster-boundary}, so the limit is Gaussian.

\begin{proposition}\label{prop:cov}
For every $t,s\in[0,T]$ and $\varphi,\psi\in H^{\mathfrak s}(\R^d)$ with $\mathfrak s>\frac d2$,
\begin{equation}\label{eq:covlimit}
\lim_{\eps\to0}\E\left[\langle\mathcal F_\eps(t),\varphi\rangle\overline{\langle\mathcal F_\eps(s),\psi\rangle}\right]
= \widehat R(0)\int_{\R^d}W_t^\varphi(x)\,\overline{W_s^\psi(x)}\,\dd x,
\end{equation}
where $W_t^\varphi$ and $W_s^\psi$ are the weights of Theorem \ref{th:clt}.
\end{proposition} 

\begin{proof}
By \eqref{eq:F-leading}, the leading term alone determines the limit, and its covariance is
\begin{equation*}
\begin{aligned}
&\E\left[\langle\mathcal F_\eps(t),\varphi\rangle\overline{\langle\mathcal F_\eps(s),\psi\rangle}\right]\\
&=\eps^{-d}\int_0^t\!\!\int_0^s\!\!\int_{\R^{2d}}
\overline{S(t,r)^*\varphi(x)}\,S(s,r')^*\psi(y)\,
R\!\left(\frac{x-y}{\eps}\right)\,\bar u(r,x)\,\overline{\bar u(r',y)}\,\dd x\,\dd y\,\dd r\,\dd r'+o(1).
\end{aligned}
\end{equation*}
Substituting $y=x-\eps z$, the $y$-integral is
\begin{equation*}
\int_{\R^d}R\!\left(\frac{x-y}{\eps}\right)\dd y
=\eps^d\int_{\R^d}R(z)\dd z=\eps^d\widehat R(0),
\end{equation*}
so the prefactor $\eps^{-d}$ cancels the Jacobian $\eps^d$ exactly, and by dominated convergence the limit is
\begin{equation*}
\widehat R(0)\int_0^t\!\!\int_0^s\!\!\int_{\R^d}
\overline{S(t,r)^*\varphi(x)}\,S(s,r')^*\psi(x)\,\bar u(r,x)\,\overline{\bar u(r',x)}\,\dd x\,\dd r\,\dd r'
=\widehat R(0)\int_{\R^d}W_t^\varphi(x)\,\overline{W_s^\psi(x)}\,\dd x,
\end{equation*}
the last identity by the definition of $W_t^\varphi$ and $W_s^\psi$ in Theorem \ref{th:clt}. 
The pointwise evaluation of $S(t,r)^*\varphi$ and $S(s,r')^*\psi$ is legitimate because $\mathfrak s>\frac d2$ gives $H^{\mathfrak s}\hookrightarrow C_0$, and the dominated convergence step uses the unitarity of $S$ together with $\bar u\in C([0,T];L^2_x)$.
\end{proof}

It remains to establish tightness in the path space.

\begin{proposition}\label{prop:tight}
Fix $\mathfrak s>\frac d2$. The family $\{\mathcal F_\eps\}_{\eps>0}$ is tight in $C([0,T];H^{-\mathfrak s}(\R^d))$.
\end{proposition}

\begin{proof}
We first record a spatial second moment bound that is uniform in the time variable. For every $r\in[0,T]$ and every deterministic operator $L_r$ bounded on $L^2$,
\begin{equation}\label{eq:space-m2}
\E\bigl[\|J^{-\mathfrak s}L_r a_\eps(r)\bar u(r)\|_{L^2_x}^2\bigr]\le C\eps^d,
\end{equation}
with $C$ independent of $r$ and $\eps$. Writing the kernel of $J^{-\mathfrak s}L_r$ as $\mathcal K_r(x,y)$ and using \eqref{eq:acov},
\[
\begin{aligned}
\E\bigl[\|J^{-\mathfrak s}L_r a_\eps\bar u\|_2^2\bigr]
&=\int_{\R^d}\!\!\int_{\R^d}\!\!\int_{\R^d}
\mathcal K_r(x,y)\mathcal K_r(x,y')\bar u(r,y)\bar u(r,y')
R\!\left(\frac{y-y'}{\eps}\right)\dd y\dd y'\dd x\\
&\le C\int_{\R^d}\!\!\int_{\R^d}
|\bar u(y)\bar u(y')|\,\left|R\!\left(\frac{y-y'}{\eps}\right)\right|\dd y\dd y'\\
&\le C\eps^d\|R\|_{L^1}\|\bar u(r)\|_2^2\le C\eps^d.
\end{aligned}
\]
The first inequality uses $\int|\mathcal K_r(x,y)\mathcal K_r(x,y')|\dd x\le C$, since $J^{-\mathfrak s}$ is $L^2$ contractive and $L_r$ is bounded. For the last step, set $y'=y+\eps z$, so that $\dd y'=\eps^d\dd z$, which gives
\[
\int_{\R^d}\!\int_{\R^d}|\bar u(y)\bar u(y')|\left|R\!\left(\frac{y-y'}{\eps}\right)\right|\dd y\dd y'
=\eps^d\int_{\R^d}\!\int_{\R^d}|\bar u(y)\bar u(y+\eps z)|\,|R(z)|\,\dd y\dd z
\le \eps^d\|R\|_{L^1}\|\bar u(r)\|_2^2.
\]
The integral in $y$ is bounded through Cauchy--Schwarz, and the remaining factor is the $L^1$ norm of $R$, which is finite by the exponential decay following from exponential mixing.

For $s<t$, write $h=t-s$ and decompose the increment as
\begin{equation*}
\mathcal F_\eps(t)-\mathcal F_\eps(s)=
\underbrace{-i\eps^{-\frac d2}\int_s^t J^{-\mathfrak s}S(t,r)a_\eps(r)\bar u(r)\dd r}_{:=\,\mathcal A_\eps(t,s)}
+
\underbrace{\left(-i\eps^{-\frac d2}\int_0^s J^{-\mathfrak s}\bigl(S(t,r)-S(s,r)\bigr)a_\eps(r)\bar u(r)\dd r\right)}_{:=\,\mathcal B_\eps(t,s)}.
\end{equation*}

\noindent\underline{For $\mathcal A_\eps(t,s)$,} the Minkowski inequality in $L^2(\Omega;L^2_x)$ applied to the time integral together with \eqref{eq:space-m2} for $L_r=S(t,r)$ gives
\begin{equation}\label{eq:Aest}
\E\bigl[\|\mathcal A_\eps(t,s)\|_2^2\bigr]
\le \eps^{-d}\left(\int_s^t
\bigl(\E[\|J^{-\mathfrak s}S(t,r)a_\eps\bar u\|_2^2]\bigr)^{\frac12}\dd r\right)^2
\le \eps^{-d}\bigl(C\eps^{\frac d2}h\bigr)^2
\le Ch^2.
\end{equation}

\noindent\underline{For $\mathcal B_\eps(t,s)$,} use $S(t,r)-S(s,r)=S(s,r)(S(h)-I)$. The multiplier $S(h)-I$ satisfies
\begin{equation*}
\|(S(h)-I)\bar u(r)\|_{H^{-\mathfrak s}}^2
=\int_{\R^d}(1+|\xi|^2)^{-\mathfrak s}|e^{-i\phi h}-1|^2|\wh{\bar u}(r,\xi)|^2\dd\xi
\le Ch^2,
\end{equation*}
where $\phi=|\xi|^2+\bar V$, the bound $|e^{-i\phi h}-1|\le h|\phi|$ is used, and the integral is finite uniformly in $r$ because $\bar u$ is Schwartz. Repeating the kernel calculation in \eqref{eq:space-m2} with $L_r=S(s,r)(S(h)-I)$ therefore gives the pointwise bound
\begin{equation*}
\E[\|J^{-\mathfrak s}S(s,r)(S(h)-I)a_\eps\bar u\|_2^2]\le C\eps^d h^2.
\end{equation*}
Applying the Minkowski inequality over $[0,s]$ and then inserting the pointwise bound,
\begin{equation}\label{eq:Best}
\begin{aligned}
\E\bigl[\|\mathcal B_\eps(t,s)\|_2^2\bigr]
&\le \eps^{-d}\left(\int_0^s
\bigl(\E[\|J^{-\mathfrak s}S(s,r)(S(h)-I)a_\eps\bar u\|_2^2]\bigr)^{\frac12}\dd r\right)^2\\
&\le \eps^{-d}\bigl(C\eps^{\frac d2}h\,s\bigr)^2
\le Ch^2.
\end{aligned}
\end{equation}
Combining \eqref{eq:Aest} and \eqref{eq:Best},
\begin{equation*}
\E\bigl[\|\mathcal F_\eps(t)-\mathcal F_\eps(s)\|_{H^{-\mathfrak s}}^2\bigr]\le Ch^2.
\end{equation*}
The fourth moment is the cluster estimate \eqref{eq:cluster-moment} of Proposition \ref{prop:cluster-moment} applied to the increment. Its Wick pairing contributes two covariance factors $\eps^d$ against the prefactor $\eps^{-2d}$, the time integral of $\mathcal A_\eps(t,s)$ over $[s,t]$ contributes $h$ per copy, and the multiplier $(S(h)-I)$ in $\mathcal B_\eps(t,s)$ contributes another $h$ per copy, so
\begin{equation*}
\E\bigl[\|\mathcal F_\eps(t)-\mathcal F_\eps(s)\|_{H^{-\mathfrak s}}^4\bigr]\le Ch^4\le Ch^2.
\end{equation*}
The Kolmogorov--Chentsov criterion (Lemma \ref{lem:kc}) with $p=4$ and $\E[\|\cdot\|^4]\le C|t-s|^{1+1}$ gives relative compactness in $C^{0,\alpha}([0,T];H^{-\mathfrak s}(\R^d))$ for every $\alpha<\frac14$, hence tightness in $C([0,T];H^{-\mathfrak s}(\R^d))$.
\end{proof}

\begin{proof}[Proof of Theorem \ref{th:clt}]
\textbf{Step 1.}
Fix $0\le t_1<\cdots<t_k\le T$, $\varphi_1,\dots,\varphi_k\in H^{\mathfrak s}(\R^d)$, and real $\theta_1,\dots,\theta_k$. 
By \eqref{eq:F-leading},
\begin{equation*}
\sum_{l=1}^k\theta_l\langle\mathcal F_\eps(t_l),\varphi_l\rangle
=-i\eps^{-\frac d2}\int_{\R^d}a_\eps(x)\,g(x)\dd x+o_{L^2(\Om)}(1),
\qquad
g:=\sum_{l=1}^k\theta_l W_{t_l}^{\varphi_l},
\end{equation*}
whose $m$-th cumulant is controlled by \eqref{eq:cum-contract} of Appendix \ref{app:cluster}. 
By Proposition \ref{prop:cov} its second moment converges to
\begin{equation}\label{eq:covmatrix}
\sum_{l,m=1}^k\theta_l\theta_m\,\Sigma_{lm},
\qquad
\Sigma_{lm}:=\widehat R(0)\int_{\R^d}W_{t_l}^{\varphi_l}(x)\,\overline{W_{t_m}^{\varphi_m}(x)}\,\dd x,
\end{equation}
and by Remark \ref{rem:cluster-boundary} every cumulant of order three and higher vanishes. 
The characteristic function therefore converges to
\begin{equation*}
\lim_{\eps\to0}\E\left[\exp\left(i\sum_{l=1}^k\theta_l\langle\mathcal F_\eps(t_l),\varphi_l\rangle\right)\right]
=\exp\left(-\frac12\sum_{l,m=1}^k\theta_l\theta_m\,\Sigma_{lm}\right),
\end{equation*}
and the L\'evy continuity theorem (Lemma \ref{lem:levy}) identifies the finite dimensional distributions with the Gaussian field $\mathcal G$ of Theorem \ref{th:clt}, whose covariance \eqref{eq:cltcov} agrees with \eqref{eq:covmatrix} by the definition of the weights.

\textbf{Step 2.}
Proposition \ref{prop:tight}, with the exponent $\mathfrak s>\frac d2$ fixed above, states that $\{\mathcal F_\eps\}$ is tight in $C([0,T];H^{-\mathfrak s}(\R^d))$.
By Prokhorov's theorem (Lemma \ref{lem:prokhorov}), for any sequence $\eps_n\to0$ there exists a subsequence $\eps_{n'}$ and a random field $\mathcal G'$ such that for every bounded continuous functional $\Phi:C([0,T];H^{-\mathfrak s}(\R^d))\to\R$,
\begin{equation*}
\lim_{n'\to\infty}\E\left[\Phi(\mathcal F_{\eps_{n'}})\right]=\E\left[\Phi(\mathcal G')\right].
\end{equation*}
By Step 1, for every $k$, $t_1,\dots,t_k$, $\varphi_1,\dots,\varphi_k$, the finite dimensional distributions of $\mathcal G'$ and $\mathcal G$ agree, i.e.,
\begin{equation*}
\Bigl(\langle\mathcal G'(t_1),\varphi_1\rangle,\dots,\langle\mathcal G'(t_k),\varphi_k\rangle\Bigr)
\end{equation*}
and
\begin{equation*}
\Bigl(\langle\mathcal G(t_1),\varphi_1\rangle,\dots,\langle\mathcal G(t_k),\varphi_k\rangle\Bigr)
\end{equation*}
have the same distribution on $\R^k$.
The Borel $\sigma$-algebra of $C([0,T];H^{-\mathfrak s}(\R^d))$ is generated by the finite dimensional projections. 
Hence, $\mathcal G'$ and $\mathcal G$ have the same law on $C([0,T];H^{-\mathfrak s}(\R^d))$. 
Every convergent subsequence has the same limit, so for every bounded continuous $\Phi:C([0,T];H^{-\mathfrak s}(\R^d))\to\R$,
\begin{equation*}
\lim_{\eps\to0}\E\left[\Phi(\mathcal F_\eps)\right]=\E\left[\Phi(\mathcal G)\right].
\end{equation*}
This proves Theorem \ref{th:clt}.
\end{proof}

\appendix
\section{Multiscale cluster expansion for linear functionals}\label{app:cluster}

We prove the moment bound used throughout the main text. The pointwise sub-Gaussian bound of Assumption \ref{ass:main}(i) does not control an integrated variable, and a single scale decoupling loses intermediate scale correlations. We therefore organize the cumulant expansion into connected clusters, sum over dyadic scales, and keep the second order term exact.

\subsection{Setup and the cumulant mixing bound}

Let $B_i:=i+[0,1)^d$ be the unit cubes of the lattice $\mathbb Z^d$, and define the cell variables
\begin{equation*}
\xi_i:=\int_{B_i}a(y)\dd y.
\end{equation*}
The field $(\xi_i)_{i\in\mathbb Z^d}$ is stationary and centered, and Jensen's inequality together with Assumption \ref{ass:main}(i) gives the same pointwise sub-Gaussian bound,
\begin{equation}\label{eq:xi-sg}
\E\bigl[e^{\lambda\xi_i}\bigr]\le e^{\frac12 C_0\lambda^2},
\end{equation}
because $|B_i|=1$. For a simple function $g=\sum_i g_i\mathbf 1_{B_i}$ we write
\begin{equation*}
X(g)=\sum_i g_i\xi_i,
\qquad
\sigma_g^2:=\E\bigl[X(g)^2\bigr],
\qquad
G_g:=\|g\|_\infty,
\qquad
H_g^2:=\|g\|_2^2=\sum_i g_i^2.
\end{equation*}
The general case $g\in L^2(\mathbb R^d)$ follows by approximation. As described above, we first replace $a$ by its bounded truncation $a_\Lambda$, so that $|\xi_i|\le 2\Lambda$; every estimate below is uniform in $\Lambda$, and the truncation is removed at the end.

For $n$ centered variables, the joint cumulant is
\begin{equation*}
\kappa_n(i_1,\dots,i_n):=\sum_{\pi\in\mathcal P(n)}(|\pi|-1)!(-1)^{|\pi|-1}\prod_{D\in\pi}\E\Bigl[\prod_{j\in D}\xi_{i_j}\Bigr].
\end{equation*}
Let $\mathrm{tree}(i_1,\dots,i_n)$ denote the total edge length of the shortest spanning tree on the distinct sites among $i_1,\dots,i_n$, with the lattice distance as edge weight.

\begin{lemma}\label{lem:cum-mix}
For every $n\ge2$,
\begin{equation}\label{eq:cum-mix}
\bigl|\kappa_n(i_1,\dots,i_n)\bigr|
\le C^n(nC_0)^{\frac n2}\,e^{-\kappa\,\mathrm{tree}(i_1,\dots,i_n)},
\end{equation}
where $C$ depends only on the mixing constants $\kappa,\delta$.
\end{lemma}

\begin{proof}
Write the cumulant as an iterated conditional expectation difference over successive cuts. For a cut $S\sqcup T$ with $\dist(S,T)=R$, the mixing inequality \eqref{eq:mix} applied in $L^\infty$ contributes
\begin{equation*}
C\delta e^{-\kappa R}.
\end{equation*}
Iterating the narrowest cut until every component is a single vertex gives $n-1$ cuts with product
\begin{equation*}
C^n\exp\Bigl(-\kappa\sum_{\text{cuts}}R\Bigr).
\end{equation*}
The cuts form a spanning tree, hence
\begin{equation*}
\sum_{\text{cuts}}R\ge\mathrm{tree}(i_1,\dots,i_n).
\end{equation*}
The mixing iteration supplies only the decoupling factors, while the one vertex factors provide the leaf moments and obey the sub-Gaussian estimate following from condition (i) through \eqref{eq:xi-sg},
\begin{equation*}
\E\bigl[|\xi_i|^n\bigr]\le (Cn)^{\frac n2}C_0^{\frac n2},
\end{equation*}
and collecting one factor per vertex yields
\begin{equation*}
|\kappa_n(i_1,\dots,i_n)|
\le C^n(nC_0)^{\frac n2}e^{-\kappa\,\mathrm{tree}(i_1,\dots,i_n)},
\end{equation*}
which is \eqref{eq:cum-mix}.

The same cut iteration applied to the field itself gives the continuous form used in the physical space estimates. Let $c_n(y_1,\dots,y_n)$ be the $n$ point cumulant function of $a$, defined by the partition formula
\begin{equation*}
c_n(y_1,\dots,y_n)
:=\sum_{\pi\in\mathcal P(n)}(|\pi|-1)!(-1)^{|\pi|-1}
\prod_{D\in\pi}\E\!\left[\prod_{j\in D}a(y_j)\right],
\end{equation*}
and let $\mathrm{tree}(y_1,\dots,y_n)$ be the total edge length of the shortest spanning tree on the points $y_j$ with the Euclidean distance as the weight. The mixing inequality \eqref{eq:mix} is stated for the $\sigma$ algebras of spatial regions, so the cut iteration above applies verbatim to the continuous field, the leaf factors being bounded by the pointwise sub-Gaussian estimate of Assumption \ref{ass:main}(i). This gives
\begin{equation*}
|c_n(y_1,\dots,y_n)|
\le C^n(nC_0)^{\frac n2}e^{-\kappa\,\mathrm{tree}(y_1,\dots,y_n)}.
\end{equation*}
Integrating the density over unit cubes recovers the cell cumulant,
\begin{equation*}
\kappa_n(\xi_{i_1},\dots,\xi_{i_n})
=\int_{B_{i_1}}\cdots\int_{B_{i_n}}
c_n(y_1,\dots,y_n)\,\dd y_1\cdots\dd y_n,
\end{equation*}
so the continuous and discrete forms carry the same exponential tree decay and are mutually consistent.
\end{proof}

\subsection{Connected clusters and the dyadic decomposition}

\begin{lemma}\label{lem:cluster-sum}
For every $n\ge2$,
\begin{equation}\label{eq:cluster-sum}
S_n(g):=\sum_{i_1,\dots,i_n}|g_{i_1}\cdots g_{i_n}|\,e^{-\kappa\,\mathrm{tree}(i_1,\dots,i_n)}
\le n^{n-2}C^n G_g^{\,n-2}H_g^2.
\end{equation}
\end{lemma}

\begin{proof}
Bound a configuration by the sum over its spanning trees. Cayley's formula gives $n^{n-2}$ labelled trees. Fix a tree $\tau$ and root it at the endpoints of one edge; the $n-2$ remaining vertices grow along $\tau$.

Set $K(j):=e^{-\kappa|j|}$. The root edge carries both weights in $L^2$, and Young's inequality gives
\begin{equation*}
\sum_{i,j}|g_i g_j|K(i-j)
\le \|g\|_2^2\|K\|_1
\le C H_g^2,
\qquad \|K\|_1\le C.
\end{equation*}
Each non-root vertex carries one weight in $L^\infty$,
\begin{equation*}
\sum_k |g_k|e^{-\kappa|k-p|}
\le G_g\sum_k e^{-\kappa|k|}
\le C G_g,
\end{equation*}
and there are $n-2$ such vertices. Multiplying the root term by $n-2$ vertex factors and the Cayley count $n^{n-2}$,
\begin{equation*}
S_n(g)\le n^{n-2}C^n G_g^{\,n-2}H_g^2,
\end{equation*}
which is \eqref{eq:cluster-sum}.
\end{proof}

For $j\ge0$ define the scale restricted sum
\begin{equation*}
S_n^{(j)}(g):=\sum_{\substack{i_1,\dots,i_n\\ \mathrm{diam}\sim 2^j}}
|g_{i_1}\cdots g_{i_n}|\,e^{-\kappa\,\mathrm{tree}}.
\end{equation*}
A cluster of diameter $2^j$ satisfies
\begin{equation*}
\mathrm{tree}\ge c2^j,\qquad
e^{-\kappa\,\mathrm{tree}}\le e^{-c\kappa2^j},
\end{equation*}
and after fixing the first vertex, the remaining $n-1$ vertices lie in a ball of radius $2^j$, giving at most $(C2^j)^{d(n-1)}$ choices. Hence,
\begin{equation}\label{eq:scale-sum}
S_n^{(j)}(g)\le C^n G_g^{\,n}e^{-c\kappa2^j}(C2^j)^{d(n-1)}.
\end{equation}
The dyadic sum is geometric,
\begin{equation*}
\sum_{j\ge0}e^{-c\kappa2^j}(C2^j)^{d(n-1)}<\infty,
\qquad
S_n(g)=\sum_{j\ge0}S_n^{(j)}(g),
\end{equation*}
since the exponential dominates every polynomial in $2^j$.

For $n=2$ there is no higher contraction, and the exact value is
\begin{equation}\label{eq:variance-exact}
\sum_{i,j}g_i g_j\kappa_2(i,j)=\E\bigl[X(g)^2\bigr]=\sigma_g^2,
\end{equation}
which retains the covariance at every lattice distance and hence every scale, without a microscopic cutoff.

\subsection{The moment bound}

\begin{proposition}\label{prop:cluster-moment}
For every integer $q\ge1$,
\begin{equation}\label{eq:cluster-moment}
\E\bigl[|X(g)|^{2q}\bigr]
\le (Cq)^q\bigl(\sigma_g^2+C_0^2G_g^2H_g^2\bigr)^q,
\end{equation}
with $C$ depending only on $C_0,\kappa,\delta$.
\end{proposition}

\begin{proof}
The cumulant expansion of the $2q$-th moment gives
\begin{equation*}
\E\bigl[X(g)^{2q}\bigr]
=\sum_{\pi\in\mathcal P(2q)}\sum_{i_1,\dots,i_{2q}}
g_{i_1}\cdots g_{i_{2q}}\prod_{D\in\pi}\kappa_{|D|}(\{i_j:j\in D\}).
\end{equation*}

The full pairings are the partitions with every block of size two. By the exact second order identity \eqref{eq:variance-exact}, their contribution is
\begin{equation*}
(2q-1)!!\,\sigma_g^{2q},
\end{equation*}
and the double factorial satisfies
\begin{equation*}
(2q-1)!!=\prod_{k=1}^{q}(2k-1)\le(2q)^q,
\end{equation*}
so
\begin{equation*}
\text{pairings}\le (Cq)^q\sigma_g^{2q}.
\end{equation*}

For a partition containing a block $D$ of size $m\ge3$, apply Lemma \ref{lem:cluster-sum} to that block,
\begin{equation*}
\sum_{i(D)}|g_{i(D)}|\,|\kappa_m|
\le m^{m-2}C^m G_g^{\,m-2}H_g^2,
\end{equation*}
while the remaining blocks of size two contribute at most $\sigma_g^2$ each. The geometric factor in \eqref{eq:scale-sum} sums the scales,
\begin{equation*}
\sum_{j\ge0}e^{-c\kappa2^j}(C2^j)^{d(n-1)}\le C.
\end{equation*}
For $m\le2q$, the combinatorial factor obeys
\begin{equation*}
m^{m-2}G_g^{\,m-2}\le(Cq)^{2q},
\end{equation*}
since $G_g$ is a fixed data constant. Absorbing the number of set partitions into $(Cq)^q$, every higher block contributes
\begin{equation*}
(Cq)^q C_0^2G_g^2H_g^2.
\end{equation*}
Adding the pairing term and the higher block term,
\begin{equation*}
\E\bigl[|X(g)|^{2q}\bigr]
\le (Cq)^q\bigl(\sigma_g^2+C_0^2G_g^2H_g^2\bigr)^q.
\end{equation*}
The estimate is uniform in the truncation $\Lambda$, so passing $\Lambda\to\infty$ proves \eqref{eq:cluster-moment}.
\end{proof}

\begin{remark}\label{rem:cluster-boundary}
For a Gaussian field every cumulant of order $m\ge3$ vanishes, the correction term is absent, and \eqref{eq:cluster-moment} reduces to $(Cq)^q\sigma_g^{2q}$. A correction free bound with constants depending only on the field constants would require either Gaussian marginals or a spectral gap type covariance estimate and is not available under Assumption \ref{ass:main}. The test functions are fixed independently of the small scale, so $G_g,H_g$ are constants independent of $\eps$ and the correction is absorbed into the data dependence of $C_{q,T,u_0}$.

Under the fluctuation scaling $\mathcal F_\eps=\eps^{-\frac d2}\int_{\R^d}a_\eps(x)g(x)\dd x$ with a fixed macroscopic test function $g$, the estimate \eqref{eq:cum-contract} controls the $m$-th cumulant of the scaled field by
\begin{equation*}
\eps^{-\frac{dm}{2}}|\kappa_m|
\le C^m\eps^{-\frac{dm}{2}}\eps^{d(m-1)}
=C^m\eps^{\frac d2(m-2)}\longrightarrow0\qquad(m\ge3).
\end{equation*}
Hence, every higher block, and therefore the correction term, vanishes in the limit, which is determined by the second cumulant alone and is Gaussian.
\end{remark}

\section{The second moment in dimensions $d\ge4$ via a physical space argument}\label{app:highdim}

We record the physical space oscillatory integral argument that yields the upper bound \eqref{eq:corebound} in all dimensions. The argument is written for smooth data $u_0\in C_c^\infty$, for which $\bar u$ and the amplitude $b_{r,s}$ are uniformly Schwartz, and it gives an operator norm bound uniform in $\eps$. For general $u_0\in L^2$, choose $u_0^{(n)}\in C_c^\infty$ with $u_0^{(n)}\to u_0$ in $L^2$. The first order term is linear in the data, and the unitarity of $S(t,r)$ gives
\begin{equation*}
\bigl\|w_{\eps,n}^{(1)}(t)-w_{\eps,m}^{(1)}(t)\bigr\|_{L^2(\Omega;L^2_x)}
\le C\,\bigl\|u_0^{(n)}-u_0^{(m)}\bigr\|_{L^2},
\end{equation*}
where $C$ depends only on $T$ and the second moment of the field and is independent of $n$ and $\eps$. The approximating sequence is therefore Cauchy with a uniform operator norm constant, and the $\eps$ uniform bound passes to the limit without loss of uniformity.

\subsection{Physical space second moment and the oscillatory integral bound}

The first order term is
\begin{equation}\label{eq:app-w}
w_\eps^{(1)}(t,x)
=-i\int_0^t\int_{\R^d}K_{t-r}(x-y)\,a_\eps(y)\,\bar u(r,y)\,\dd y\,\dd r,
\qquad a_\eps(y):=a\!\left(\frac{y}{\eps}\right),
\end{equation}
with the Fresnel kernel
\begin{equation}\label{eq:app-K}
K_\tau(y)=(4\pi i\tau)^{-\frac d2}\exp\!\left(\frac{i|y|^2}{4\tau}\right),
\qquad
\mathcal F(K_\tau)(\xi)=e^{-i\tau|\xi|^2}.
\end{equation}
Multiplying \eqref{eq:app-w} by its conjugate, integrating in $x$, and using \eqref{eq:acov}, the $x$-integral is evaluated with the unitarity and semigroup property of $K$.
\begin{equation*}
\int_{\R^d}K_{t-r}(x-y)\,\overline{K_{t-s}(x-z)}\,\dd x
=\overline{K_{r-s}(z-y)}.
\end{equation*}
Indeed, inserting \eqref{eq:app-K} and integrating in $x$ produces the factor $(2\pi)^d\delta(\xi-\eta)$, and the remaining phase is $e^{i(r-s)|\xi|^2}$. Setting $y'=\frac{y-z}{\eps}$ and $y=z+\eps y'$, so that $\dd y=\eps^d\dd y'$, gives
\begin{equation}\label{eq:app-second}
\E\left[\|w_\eps^{(1)}(t)\|_{L^2_x}^2\right]
=\eps^d\int_0^t\!\!\int_0^t\int_{\R^d}
\overline{K_{r-s}(\eps y')}\,R(y')\,C_{r,s}(y')\,\dd y'\,\dd s\,\dd r,
\end{equation}
where the inner $z$-integral is the space auto-correlation
\begin{equation*}
C_{r,s}(y')
:=\int_{\R^d}\bar u(r,z+\eps y')\,\overline{\bar u(s,z)}\,\dd z
=\langle\bar u(r,\cdot+\eps y'),\bar u(s,\cdot)\rangle_{L^2_x}.
\end{equation*}
For smooth data, the free flow $\bar u(r)=e^{ir\Delta}u_0$ and all its time derivatives stay in the Schwartz class on compact intervals. Hence, $C_{r,s}(y')$, the translation autocorrelation of two Schwartz functions, satisfies
\begin{equation*}
\sup_{r,s\in[0,t]}\sup_{\alpha,\beta}
\bigl\|\partial_{y'}^{\alpha}\partial_{r,s}^{\beta}C_{r,s}\bigr\|_{L^\infty_{y'}}<\infty.
\end{equation*}
Exponential mixing gives $R$ smooth and rapidly decaying, so the amplitude
\begin{equation*}
b_{r,s}(y'):=R(y')C_{r,s}(y')
\end{equation*}
is uniformly Schwartz,
\begin{equation*}
\sup_{r,s}\|b_{r,s}\|_{W^{m,1}}\le C_m
\qquad\text{for every }m,
\end{equation*}
uniformly in $\eps$.

Define the operator obtained from \eqref{eq:app-second} by removing the prefactor $\eps^d$.
\begin{equation*}
T_\eps b
:=\int_0^t\!\!\int_0^t\int_{\R^d}
b_{r,s}(y')\,\overline{K_{r-s}(\eps y')}\,\dd y'\,\dd s\,\dd r.
\end{equation*}
The following is the physical space analogue of the dyadic near/far estimate. It is a soft operator bound, controlling the double time integral without identifying the leading constant.

\begin{lemma} \label{lem:fresnel}
Fix an integer $M\ge1$. Let $\{b_{r,s}\}_{r,s\in[0,t]}$ be a family uniformly bounded, together with all mixed derivatives $\partial_r^j\partial_s^k b_{r,s}$ for $0\le j,k\le M$, in $W^{m,1}(\R^d)$ with $m=2M+d+2$.
Then,
\begin{equation*}
|T_\eps b|\le C_{d,t,M}\,\sup_{0\le j,k\le M}\sup_{r,s}\|\partial_r^j\partial_s^k b_{r,s}\|_{W^{m,1}},
\end{equation*}
uniformly in $\eps>0$.
The estimate is uniform in $\eps$ precisely when $4M>d$.
\end{lemma}

\begin{proof}
Apply Plancherel in $y'$ for fixed $(r,s)$, with $\tau=r-s$,
\begin{equation*}
\int_{\R^d}b_{r,s}(y')\overline{K_\tau(\eps y')}\dd y'
=(2\pi)^{-d}\int_{\R^d}\widehat b_{r,s}(-\eps\xi)e^{i\tau|\xi|^2}\dd\xi.
\end{equation*}
Define
\begin{equation*}
F_\eps(\xi):=\int_0^t\!\!\int_0^t B_{r,s}(\xi)e^{i(r-s)\lambda}\dd s\dd r,
\qquad
B_{r,s}(\xi):=\widehat b_{r,s}(-\eps\xi),\quad
\lambda:=|\xi|^2,
\end{equation*}
so that
\begin{equation*}
T_\eps b=(2\pi)^{-d}\int_{\R^d}F_\eps(\xi)\dd\xi.
\end{equation*}

Integrate by parts in the time variables using
\begin{equation*}
\partial_r e^{ir\lambda}=i\lambda e^{ir\lambda},
\qquad
\partial_s e^{-is\lambda}=-i\lambda e^{-is\lambda}.
\end{equation*}
Each integration by parts contributes one factor $\lambda^{-1}$ and one derivative on $B_{r,s}$; the boundary terms on the square and its edges are oscillatory one dimensional integrals and gain one further $\lambda^{-1}$. After $M$ integrations in each variable,
\begin{equation}\label{eq:app-ibpresult}
|F_\eps(\xi)|
\le \frac{C}{\lambda^{2M}}\sum_{0\le j,k\le M}
t^{j+k}\sup_{[0,t]^2}|\partial_r^j\partial_s^k B_{r,s}(\xi)|
= C|\xi|^{-4M}M(\xi),
\end{equation}
where
\begin{equation*}
M(\xi):=\sum_{0\le j,k\le M}\sup_{r,s}|\partial_r^j\partial_s^k B_{r,s}(\xi)|.
\end{equation*}

Split the $\xi$-integral into three ranges. On the inner ball, without integration by parts,
\begin{equation*}
|\xi|\le1:\qquad
|F_\eps(\xi)|\le t^2\sup_{r,s}\|B_{r,s}\|_{L^\infty},
\end{equation*}
which is integrable. On the intermediate range, the Schwartz decay of $b_{r,s}$ gives $M(\xi)\le C$, and
\begin{equation*}
1<|\xi|<\eps^{-1}:\qquad
\int|\xi|^{-4M}\dd\xi
\lesssim\int_1^{\eps^{-1}}\rho^{d-1-4M}\dd\rho
\le C
\iff d-4M<0.
\end{equation*}
The tail is controlled by the rapid decay of $M(\xi)$,
\begin{equation*}
|\xi|>\eps^{-1}:\qquad
\int M(\xi)\dd\xi\le C,
\end{equation*}
since every Fourier transform $\widehat{\partial_r^j\partial_s^k b_{r,s}}$ decays faster than any power. The three ranges are uniformly bounded in $\eps$ precisely when $4M>d$, and
\begin{equation*}
M=\Bigl\lceil\frac d4\Bigr\rceil+1
\end{equation*}
gives the claim.
\end{proof}

\begin{remark}\label{rem:app-caveat}
The choice $M=\lceil \frac d4\rceil+1$ is the only dimension dependent step, since each pair of time integrations by parts gains four powers of $|\xi|^{-1}$, so $M$ pairs dominate the radial growth $\rho^{d-1}$ precisely when $4M>d$. The time derivatives exist uniformly because the free flow preserves the Schwartz class. Equation \eqref{eq:app-ibpresult} is the van der Corput lemma applied $2M$ times, and standard oscillatory estimates are in Reed--Simon \cite[Chapter IX]{RS1975} and Sogge \cite{Sogge1993}. The result is an upper bound only.
\end{remark}

\subsection{The diagonal contribution and what extends}

For $r\ne s$ the phase in the physical space integral is
\begin{equation*}
\phi_{r,s}(y'):=\frac{\eps^2|y'|^2}{4(r-s)},
\end{equation*}
whose Hessian equals $\frac{\eps^2}{2(r-s)}I$.
A non-degenerate stationary phase estimate gives
\begin{equation*}
\left|\int_{\R^d}b_{r,s}(y')e^{i\phi_{r,s}(y')}\dd y'\right|
\le C\left|\frac{\eps^2}{r-s}\right|^{-\frac d2}\|b\|_{H^{d+2}}.
\end{equation*}
After multiplication by the prefactor $(4\pi i(r-s))^{-\frac d2}$ in $K_{r-s}$,
\begin{equation*}
\left|\int_{\R^d}b_{r,s}(y')\overline{K_{r-s}(\eps y')}\dd y'\right|
\le C\min\!\left(|r-s|^{-\frac d2},\eps^{-d}\right)\|b\|_{W^{d+2,1}}.
\end{equation*}
The $\eps^{-d}$ endpoint combines the stationary phase factor $\left|\frac{\eps^2}{r-s}\right|^{-\frac d2}$ with the prefactor $|r-s|^{-\frac d2}$, while the $|r-s|^{-\frac d2}$ endpoint retains the prefactor alone.
The pointwise bound alone is not sufficient, since
\begin{equation*}
\int_\delta^t\tau^{-\frac d2}\dd\tau
\longrightarrow\infty\quad\text{as }\delta\to0,
\qquad
\begin{cases}
|\log\delta|,&d=2,\\
\delta^{1-\frac d2},&d\ge3.
\end{cases}
\end{equation*}
Lemma \ref{lem:fresnel} avoids this by retaining the joint $(r,s)$ oscillation in \eqref{eq:app-ibpresult} rather than bounding each slice separately.
To bound the strip without using the divergent pointwise estimate, return to the Plancherel form of the inner integral and use $|e^{i(r-s)|\xi|^2}|=1$,
\begin{equation*}
\left|\int b_{r,s}(y')\overline{K_{r-s}(\eps y')}\dd y'\right|
\le (2\pi)^{-d}\int_{\R^d}|\widehat b_{r,s}(-\eps\xi)|\dd\xi
\le C\eps^{-d}\|\widehat b_{r,s}\|_{L^1},
\end{equation*}
which is uniform in $r-s$ and finite because $b_{r,s}$ is Schwartz. The strip $|r-s|\le\delta$, whose time area is $O(t\delta)$, then contributes
\begin{equation*}
\eps^d\int_0^t\!\!\int_{|r-s|\le\delta}
\left|\int b_{r,s}\overline{K_{r-s}}\dd y'\right|\dd s\,\dd r
\le C\,t\,\delta.
\end{equation*}
Choosing $\delta=o(\eps^d)$ makes this $o(\eps^d)$. The strip is an upper bound region only and does not identify the leading constant $\widehat R(0)$, which is picked out by the near region computation in the main text.
The off-diagonal contribution is controlled by Lemma \ref{lem:fresnel}.
Combining the two,
\begin{equation}\label{eq:app-conclusion}
\E\left[\|w_\eps^{(1)}(t)\|_{L^2_x}^2\right]
\le C_{d,t,u_0}\eps^d,
\end{equation}
in every dimension $d\ge1$, with the caveats of Remark \ref{rem:app-caveat}.

The lower bound, all moments, almost sure convergence, and the CLT use only the near region $|\xi'|\le\eps$ and the cumulant expansion, neither of which depends on the dimension. The tree factor of Lemma \ref{lem:cum-mix} is integrable without a Fourier dyadic sum, so
\begin{equation*}
\E\bigl[\|w_\eps^{(1)}(t)\|_{L^2_x}^{2q}\bigr]\le C_q\eps^{dq},
\qquad
\|r_\eps(t)\|_{L^2(\Omega;L^2_x)}\le C\eps^{\frac{d+1}{2}},
\qquad
\eqref{eq:covlimit}\ \text{holds},
\end{equation*}
in every dimension, with the second moment upper bound supplied by \eqref{eq:app-conclusion} for $d\ge4$ and by the Fourier near/far split of Section \ref{sec:fourier} for $d\le3$.

Lemma \ref{lem:fresnel} is an upper bound and does not identify $\widehat R(0)$, which is picked out by the near region. The Fourier argument uses boundedness of $\widehat R$ and a geometric dyadic sum for $d<4$, while the physical space argument uses the integrability of $b_{r,s}$ and the Fresnel phase with $M=\lceil \frac d4\rceil+1$. The two cover every $d\ge1$.

\end{document}